\documentclass[a4paper,reqno]{amsart}

\usepackage[all]{xy}
\usepackage{xspace}
\usepackage{amsmath}
\usepackage{amstext}
\usepackage{amsfonts}
\usepackage[mathscr]{euscript}
\usepackage{amscd}
\usepackage{latexsym}
\usepackage{amssymb}
\usepackage{diagrams}

\newtheorem{theorem}{Theorem}[section]
\newtheorem{lemma}[theorem]{Lemma}
\newtheorem{proposition}[theorem]{Proposition}
\newtheorem{definition}[theorem]{Definition}
\newtheorem{corollary}[theorem]{Corollary}
\newtheorem{remark}[theorem]{Remark}
\theoremstyle{definition}

\newcommand{\bsy}[1]{\boldsymbol{#1}}
\newcommand{\mi}{\text{-}}
\newcommand{\enpr}{\hfill $\Box $}
\newcommand{\supar}[1]{\overset{#1}\rightarrow}

\newcommand{\prf}{\noindent\textbf{Proof. }}
\newcommand{\rw}{\rightarrow}

\newcommand{\id}{\mathrm{id}}

\newcommand{\hm}{\mathrm{Hom}}

\newcommand{\hund}[1]{{\hm}_{#1}}

\newcommand{\tiund}[1]{{\times}_{#1}}

\newcommand{\catu}{\mathrm{cat}^2\mi\mathrm{group}}

\newcommand{\ovl}[1]{\overline{#1}}

\newcommand{\pr}{\mathrm{pr}}

\newcommand{\dcl}{\mathcal{D}}

\newcommand{\ccl}{\mathcal{C}}

\newcommand{\gcl}{\mathcal{G}}

\newcommand{\ncl}{\mathcal{N}}

\newcommand{\tcl}{\mathcal{T}}

\newcommand{\hcl}{\mathcal{H}}
\newcommand{\kcl}{\mathcal{K}}
\newcommand{\pcl}{\mathcal{P}}

\newcommand{\scl}{\mathcal{S}}

\newcommand{\rt}{\rtimes}

\newcommand{\set}{\mathbf{Set}}

\newcommand{\ner}{Ner\,}
\newcommand{\cop}{\ccl^{\Delta^{op}}}
\newcommand{\cat}{\mathrm{Cat}\,}

\newcommand{\ds}{disc\,}
\newcommand{\gp}{\mathrm{Gp}}
\newcommand{\cgp}{\cat(\gp)}
\newcommand{\ccgp}{\cat^2(\gp)}
\newcommand{\pt}{\partial}

\newcommand{\tam}{Tamsamani's }
\newcommand{\dsp}{\displaystyle}
\newcommand{\pp}{\bsy{\Pi}}
\newcommand{\sext}{{\scl}/{\sim^{ext}}}
\newcommand{\kext}{{\kcl}/{\sim^{ext}}}
\newcommand{\hext}{{\hcl}/{\sim^{ext}}}
\newcommand{\ssext}{\frac{\scl}{\sim^{ext}}}
\newcommand{\kkext}{\frac{\kcl}{\sim^{ext}}}

\newcommand{\gra}{\textbf{Gray\,}}
\newcommand{\dop}{\Delta^{op}}
\newcommand{\ddop}{\Delta^{2^{op}}}
\newcommand{\bic}{Bicat\,}

\newcommand{\bi}{{Bic}\,}
\newarrow{Eq}{}{=}{=}{=}{}
\newarrow{LR}{<}{-}{-}{-}{>}
\newarrow{Dts}{}{.}{.}{.}{}
\newarrow{Embedline}{>}{-}{-}{-}{-}
\newarrow{Intoar}{C}{-}{-}{-}{>}

\newcommand{\hgra}{\hcl o(Gray\mi gpd_0)}
\newcommand{\hcon}{\hcl o(Top^{(3)}_*)}
\newcommand{\twc}{\emph{2\mi cat}}

\begin{document}

\title[Semistrict models of connected 3-types and Tamsamani's...]{Semistrict models of connected 3-types and Tamsamani's weak 3-groupoids}

\author{Simona Paoli}

\date{9 July 2006}

\address{Department of Mathematics\\ Macquarie University\\ NSW 2109 \\ Australia}

\email{simonap@ics.mq.edu.au}



\begin{abstract}
Homotopy 3-types can be modelled algebraically by Tamsamani's weak 3-groupoids as well
as, in the path connected case, by cat$^2$-groups. This paper gives a comparison between
the two models in the path-connected case. This leads to two different  semistrict
algebraic models of connected 3-types using Tamsamani's model. Both are then related to
Gray groupoids.
\end{abstract}

\maketitle

%%%%%%%%%%%%%%%%%%%%%%%%%%%%%%%%%%%%%%%%%%%%%%%%%%%%%%%%%%%%%%%%%%%%%%%%%%%%%%%%%%%%
\section{Introduction}\label{sec1}
The problem of modelling homotopy types is relevant to both homotopy theory and higher
category theory. In homotopy theory, various models exist for the path-connected case:
the cat$^n$-group model, introduced by Loday \cite{lod} and later developed in
\cite{bcd},\cite{por} generalized the earlier work of Whitehead on crossed modules
\cite{whit}; another model was built by Carrasco and Cegarra \cite{cg}.

In higher category theory, homotopy models serve as a ``test" for a good definition of
weak higher category, which should give a model of $n$-types in the weak $n$-groupoid
case. This property has been proved in \cite{tam},\cite{tam1} for the \tam model of weak
$n$-categories.

Tamsamani's model and cat$^n$-groups are both multi-simplicial models but they have
distinctly different features: the first is a strict but cubical higher categorical
structure, the second is a weak and globular one. Yet they both encode (path-connected)
$n$-types. Their differences as well as similarities stimulated our interest in searching
for a direct comparison between the two. This paper solves the problem for the case
$n=3$; that is, for connected 3-types.

As a result of this comparison we find a model of connected 3-types through a subcategory
$\mathcal{H}$ of Tamsamani's weak 3-groupoids whose objects are not, in general, strict,
but which are `less weak' than the ones used by Tamsamani to model connected 3-types
(Definition \ref{sec4.def1}). Thus objects of $\mathcal{H}$ are `semistrict' structures.
Yet the subcategory $\mathcal{H}$ is not isomorphic to the category of Gray groupoids
with one object (Remark \ref{sec5.rem1}). Thus our comparison yields a new semistrict
model of connected 3-types and at the same time gives a refinement of Tamsamani's result
in this case (Theorem \ref{sec5.the1} and Corollary \ref{sec5.cor1}).

In a generic Tamsamani's weak 3-groupoid, there are two fixed directions in the
3-simplicial set in which the Segal maps are equivalences rather than isomorphisms. In
the semistrict model $\mathcal{H}$ the Segal maps in one fixed direction between the two
are isomorphisms. Choosing the other direction for the isomorphisms of the Segal maps
leads to a different semistrict subcategory $\mathcal{K}$ (Definition \ref{sec5a.def2}).
It is then natural to ask whether $\mathcal{K}$ constitutes another semistrict model of
connected 3-types. The answer is positive (Theorem \ref{sec5a.the1}). Finally, using
simultaneously the semistrict models $\mathcal{H}$ and $\mathcal{K}$, the connection with
Gray groupoids (with one object) is established (Theorem \ref{sec6.the1}).

The paper is organized as follows.  In Section \ref{sec1a} we give a review of the
definitions and main properties of the \tam model and of the $\catu$ model which we are
going to use. In Section \ref{sec2} we prove a general fact about internal categories in
categories of groups with operations. In Section \ref{sec3} we specialize this to the
case where $\ccl$ is the category $\cat(\gp)$ of internal categories in groups. This
yields Proposition \ref{sec3.pro1}, which is the key to the passage from a strict cubical
structure to a weak globular one. Proposition \ref{sec3.pro1} says that we can represent
a connected 3-type, up to homotopy, through a cat$^2$-group with the property that, as an
internal category in $\cat(\gp)$, its object of objects is projective in $\cat(\gp)$.

The form of the projective objects in $\cat(\gp)$, known in the literature, allows one to
associate to a $\catu$ $\gcl$ of the type above a bisimplicial group $ds\,\ncl\gcl$,
which we call the \emph{discrete multinerve} of $\gcl$. The discrete multinerve of $\gcl$
has the following properties (Lemma \ref{sec3.lem3}).
\begin{itemize}
  \item [i)] For each $n\geq 0$, $(ds\,\ncl\gcl)_n:\dop\rw\gp$ is the nerve of an object
  of $\cat(\gp)$.
  \item [ii)] $(ds\,\ncl\gcl)_0$ is a constant functor.
  \item [iii)] The Segal maps are weak equivalences of simplicial groups.
\end{itemize}
Further, the classifying spaces of $\gcl$ and of its discrete multinerve are weakly
homotopy equivalent.

In general, we call a bisimplicial group with properties i), ii), iii), an \emph{internal
2-nerve} (Definition \ref{sec3.def1}). We can think of an internal 2-nerve as a version
of \tam weak 2-nerves internalized in the category of Groups. Proposition \ref{sec3.pro1}
allows use of the discrete multinerve to define the functor \emph{discretization}
\begin{equation*}
  \ds:\ccgp/\!\sim\;\rw\;\dcl/\!\sim
\end{equation*}
from the localization of $\catu$s with respect to the weak equivalences to the
localization of internal 2-nerves with respect to weak equivalences. As summarized in
Proposition \ref{sec3.pro2}, the discretization functor preserves the homotopy type.

In Section \ref{sec3} we realize the passage from $\catu$ to the semistrict subcategory
$\hcl$ of \tam weak 3-groupoids (Theorem \ref{sec5.the1}). As explained in the proof of
Theorem \ref{sec5.the1}, this is achieved by composing the functor $\ds$ with the functor
$\ovl{N}:\dcl/\sim\;\;\rw\;\;\hext$ induced by the nerve functor $\gp\rw\set$.

In this proof, and later in the paper, we use the fact that a morphism of weak
3-groupoids is an external equivalence if and only if it induces a weak homotopy
equivalence of classifying spaces; in one direction, this fact was proved in
\cite[Proposition 11.2]{tam}. However, the opposite direction of it does not seem to be
stated explicitly in the literature. In the Appendix we have provided a proof of this
fact for the subcategory $\scl$ of \tam weak 3-groupoids which we use in this paper.

Section \ref{sec5a} considers a different semistrictification of \tam weak 3-groupoids
for the path-connected case, through the subcategory $\kcl$ (Definition \ref{sec5a.def2}
and Theorem \ref{sec5a.the1}). Our method uses a strictification functor from \tam weak
2-groupoids to \tam strict 2-groupoids (Definition \ref{sec5a.def1}). Section \ref{sec5a}
is independent on the comparison with cat$^2$-groups. However, in order to relate $\kcl$
to Gray groupoids, the subcategory $\hcl$ is used (Lemma \ref{sec6.lem1}). We then
associate to objects of $\hcl$ and $\kcl$, up to homotopy, a Gray groupoid (with one
object) having the same homotopy type (Theorem \ref{sec6.the1}).

\textbf{Acknowledgements} I am grateful to Michael Batanin and Clemens Berger for helpful
conversations. I also thank the members of the Australian Category Seminar for support
and encouragement. This work is supported by an Australian Research Council Postdoctoral
Fellowship (project no. DP0558598).

%%%%%%%%%%%%%%%%%%%%%%%%%%%%%%%%%%%%%%%%%%%%%%%%%%%%%%%

\section{Preliminaries.}\label{sec1a}
In this section we review the main definitions and properties of the Tamsamani model and
of the $\catu$ model. Useful references for these topics are \cite{brr}, \cite{bcd},
\cite{ccg}, \cite{lod}, \cite{por}, \cite{tam}, \cite{tam1}.

Let $\Delta$ be the simplicial category; we denote by $\mathcal{C}^{\Delta^{op}}$ the
category of simplicial objects in $\mathcal{C}$, and by $\Delta^{n^{op}}$ the product of
$n$ copies of $\Delta^{op}$ so that $\mathcal{C}^{\Delta^{n^{op}}}$ is the category of
$n$-simplicial objects in $\mathcal{C}$.

 A \emph{1-nerve} is a simplicial set which is the nerve of a small category. The
category of 1-nerves is isomorphic to Cat. A \emph{2-nerve} is a bisimplicial set
$\phi:\ddop\rw\set$, such that each $\phi_{n*}$ is a 1-nerve, $\phi_{0*}$ is constant and
the Segal maps
$\phi_{n*}\rw\phi_{1*}\tiund{\phi_{0*}}\overset{n}\cdots\tiund{\phi_{0*}}\phi_{1*}$ are
equivalences of categories. A morphism of 2-nerves is a morphism of the underlying
bisimplicial sets. We denote by $\ncl_2$ the category of \tam 2-nerves.

For a 2-nerve $\phi$, let $T\phi_{n*}$ be the set of isomorphisms classes of objects of
the category corresponding to $\phi_{n*}$. Then the simplicial set
$T\phi:\Delta^{op}\rw\set$, $T\phi([n])=T\phi_{n*}$ is a 1-nerve.

A 2-nerve $\phi$ is a \emph{weak 2-groupoid} when each $\phi_{n*}$ and $T\phi$ are nerves
of groupoids. We denote by $\tcl_2$ the category of the weak 2-groupoids. A 2-nerve is a
\emph{strict 2-groupoid} if it is a weak 2-groupoid and all Segal maps are isomorphisms.
We denote by $\tcl_2^{st}$ the category of \tam strict 2-groupoids.

Let $2\mi{gpd}$ be the category of 2-groupoids in the ordinary sense, that is the
subcategory of the category $\twc$ of those 2-categories whose cells of positive
dimension are invertible. Hence $2\mi{gpd}$ is the category of groupoid objects in
groupoids whose object of objects is discrete. By taking the nerve, we therefore have an
isomorphism
\begin{equation*}
  \nu:2\mi{gpd}\rw\tcl_2^{st}.
\end{equation*}
A morphism $f:\phi\rw\phi'$ of 2-nerves is an \emph{external equivalence} if for all
$x,y\in\phi_{0*}$ the maps $f_{(x,y)}: \phi_{(x,y)}\rw\phi'_{(fx,fy)}$ and $Tf$ are
equivalences of categories. Recall that $\phi_{(x,y)}$ are ``hom-categories" and
\begin{equation*}
  \underset{x,y\in\phi_{0*}}{\coprod}\phi_{(x,y)}=\phi_{1*}
\end{equation*}
A \emph{3-nerve} is a 3-simplicial set $\psi:\Delta^{3^{op}}\rw\set$ such that each
$\psi( [n],\mi,\mi)$ is a 2-nerve, $\psi( [0],\mi,\mi)$ is constant and the Segal maps
are external equivalences of 2-nerves. A morphism of 3-nerves is a morphism of the
underlying 3-simplicial sets.

For a 3-nerve $\psi$, the bisimplicial set $T\psi$ is a 2-nerve so we define
$T^2\psi=T(T\psi)$, which is a 1-nerve. A \emph{weak 3-groupoid} is a 3-nerve $\psi$ such
that each $\psi( [n],\mi,\mi)$ is a weak 2-groupoid and $T^2\psi$ is a groupoid.

A morphism $g:\psi\rw\psi'$ of weak 3-groupoids is an \emph{external equivalence} if for
each $x,y\in\psi( [0],\mi,\mi)$ the map $g_{(x,y)}:\psi_{(x,y)}\rw\psi'_{(fx,fy)}$ is an
external equivalence of weak 2-groupoids and if $T^2\psi$ is an equivalence of
categories. We denote by $\tcl_3$ the category of \tam weak 3-groupoids. A 3-nerve is a
\emph{strict 3-groupoid} if it is a weak 3-groupoid and all Segal maps are isomorphisms.

Tamsamani showed (\cite[Theorem 8.0]{tam}) that the homotopy category of 3-types is
equivalent to the localization of $\tcl_3$ with respect to external equivalences. Recall
that a $3$-type (resp. connected $3$-type) is a topological space (resp. connected
topological space) with trivial homotopy groups in dimension $i>3$; we denote by
$Top^{(3)}_*$ the category of connected 3-types.

Since we are dealing in this paper with the path-connected case, we restrict our
attention to the full subcategory $\scl$ of $\tcl_3$ consisting of those weak 3-groupoids
$\psi:\Delta^{3^{op}}\rw\set$ satisfying the condition that the constant functor
$\psi([0],\mi,\mi):\dop\rw\set$ takes values in the one-element set.

In fact, from the construction of the fundamental 3-groupoid functor
$\pp_3:\mathrm{Top}\rw\tcl_3$ given in \cite[Theorem 6.4]{tam} it is immediate that, when
restricting to the subcategory $\mathrm{Top}_*\subset\mathrm{Top}$ of path-connected
topological spaces, one obtains a functor $\pp_3:\mathrm{Top}_*\rw\scl$. Hence, in the
path-connected case, \tam result becomes
\begin{theorem}\rm{\cite{tam}}\em\label{sec1a.the1}
There is an equivalence of categories
\begin{equation*}
  \sext\simeq\hcon
\end{equation*}
\end{theorem}
In Section \ref{sec6} we need some results relating \tam 2-nerves and bicategories, in
order to construct a strictification functor from $\tcl_2$ to $\tcl_2^{st}$. In the
following theorem, $N\,\bic$ denotes the 2-category of bicategories, normal homomorphisms
and oplax natural transformations with identity components (see \cite{pao})
\begin{theorem}\rm{\cite[Th. 7.2]{pao}}\em\label{sec1a.the2}
    There is a fully faithful 2-functor $N\!:N\bic\rw\ncl_2$ with a left 2-adjoint $\bi$.
    The counit $\bi N\rw 1$ is invertible and each component $u:X\rw N\, \bi X$ of the
    unit is a pointwise equivalence, and $u_0,u_1$ are identities.
\end{theorem}
We now turn to the $\catu$ model. We begin by recalling some general notions.

Given a category $\ccl$ with pullbacks, an \emph{internal category in $\ccl$} consists of
a diagram
\begin{equation*}
  \phi:C_1\tiund{C_0}C_1\rTo^c C_1\pile{\rTo^{d_0}\\ \rTo^{d_1}\\ \lTo^i} C_0
\end{equation*}
where $d_0,d_1$ are ``source" and ``target" maps, $i$ is the ``identity" map, $c$ is the
``composition" map. These data satisfy the axioms of  a category; that is, the following
identities hold, where $\pi_0,\pi_1:C_1\tiund{C_0}C_1\rw C_1$ are the two projections:
\begin{align*}
   & d_{0}i=1_{C_0}=d_1 i, \qquad d_1\pi_1=d_1 c, \qquad d_0\pi_0=d_0 c\\
    &c\Bigg(\begin{matrix}1_{C_1}\\i d_0\end{matrix}\Bigg)=1_{C_1}=\Bigg(\begin{matrix} i
    d_1\\1_{C_1}\end{matrix}\Bigg), \qquad c(1_{C_1}\tiund{C_0}c)=c(c\tiund{C_0}1_{C_1}).
\end{align*}
Given  internal categories $\phi$ and $\phi'$, an \emph{internal functor}
$F:\phi\rw\phi'$ consists of a pair of morphisms $F_0:C_0\rw C'_0$, $F_1:C_1\rw C'_1$
satisfying the conditions:
\begin{align*}
   & d_0 F_1=F_0 d_0,\qquad d_1 F_1=F_0 d_1 \\
   & F_1 i=i F_0, \qquad \quad F_1 c=c(F_1\tiund{F_0}F_1).
\end{align*}
Let $\cat\ccl$ be the category of internal categories in $\ccl$ and internal functors.

In this paper we use the category $\cat\ccl$ in the case where $\ccl$ is a category of
groups with operations in the sense of \cite{orz}, \cite{por1}. Recall that this consists
of a category of groups with a set of additional operations $\Omega=\Omega_0 \cup
\Omega_1 \cup\Omega_2$, where $\Omega_i$ is the set of i-ary operations in $\Omega_i$
such that the group operations of identity, inverse and multiplication (denotes 0,+,-)
are elements of $\Omega_0$, $\Omega_1$, $\Omega_2$ respectively; it is $\Omega_0=\{0\}$
and certain compatibility conditions hold (see \cite{orz}, \cite{por1}). Further, there
is a set of identities which includes the groups laws.

When $\ccl$ is a category of groups with operations, every object of $\cat\ccl$ is an
internal groupoid, as every arrow is invertible.

An example of a category of groups with operations which will be important for us is the
category $\cat(\gp)$ of internal categories in the category of groups. The fact that this
forms a category of groups with operations follows easily from its equivalence with
cat$^1$-groups. In fact, recall \cite{lod} that a \emph{cat$^1$-group} consists of a
group $G$ with two endomorphisms $d, t:G\rw G$ satisfying
\begin{equation}\label{sec1a.eq1}
  dt=t,\quad td=d,\quad [\ker d,\ker t]=1.
\end{equation}
A morphism of cat$^1$-groups $(G,d,t)\rw(G',d',t')$ is a group homomorphism $f:G\rw G'$
such that $fd=d' f$, $ft=t' f$. The equivalence of $\cat(\gp)$ with cat$^1$-groups is
realized by associating to an object of $\cat(\gp)$
\begin{equation*}
  C_1\tiund{C_0}C_1\rTo^c C_1\pile{\rTo^{d_0}\\ \rTo^{d_1}\\ \lTo^i} C_0
\end{equation*}
the cat$^1$-group $(C_1,i\, d_0,i\,d_1)$ (see \cite{lod}).

The category of cat$^1$-groups is a category of groups with operations by setting
$\Omega_0=\{0\}$, $\Omega_1=\{-\}\cup\{t,d\}$, $\Omega_2=\{+\}$ and requiring the set of
identities to comprise the group laws and the identities (\ref{sec1a.eq1}). An
alternative proof that $\cat(\gp)$ is a category of groups with operations is found in
\cite{ccg}, using the language of crossed modules. In particular, \cite{ccg} contains a
characterization of objects in the algebraic category $\cat(\gp)$ which are projective
with respect to the class of regular epimorphisms. This characterization will be
important in Section \ref{sec3} and will be recalled explicitly there.

When $\ccl=\cat(\gp)$ the category $\cat\ccl$ is the category $\cat(\cat(\gp))$ of double
categories internal to Groups. This is equivalent to the category of $\catu$s, originally
introduced by Loday \cite{lod}.

Recall that a \emph{cat$^2$-group} consists of a group $G$ together with four
endomorphisms $t_i,d_i:G\rw G\;\;i=1,2$ such that, for all $1 \leq i,j\leq 2$
\begin{align*}
   & d_it_i=t_i,\quad t_id_i=d_i \\
   & d_it_j=t_jd_i,\quad d_id_j=d_jd_i,\quad t_it_j=t_jt_i,\quad i\neq j\\
   & [\ker d_i,\ker t_i]=1
\end{align*}
A \emph{morphism} of $\catu$s $(G,d_i,t_i)\rw(G',d'_i,t'_i)$ consists of a group
homomorphism $f:G\rw G'$ such that $fd_i=d'_i f$, $ft_i=t'_if$, $1\leq i\leq 2$. It is
well known that the category of $\catu$s is isomorphic to the category $\cat(\cat(\gp))$.
Because of this, when in this paper we talk about the category of $\catu$s, we always
mean the category $\cat(\cat(\gp))$, which we denote by $\ccgp$.

Given an object $\gcl$ of $\cat(\cat(\gp))$, by applying the nerve functor twice one
obtains a bisimplicial group $\ncl\gcl$, called the \emph{multinerve} of $\gcl$. The
\emph{classifying space} $B\gcl$ of the $\catu$ $\gcl$ is, by definition, the classifying
space of its multinerve. It can be shown \cite{bcd} that $B\gcl$ is a connected 3-type.

A morphism of $\catu$s is a \emph{weak equivalence} if it induces a weak homotopy
equivalence of classifying spaces.

In \cite{bcd} a \emph{fundamental cat$^2$-group} functor $\pcl:\mathrm{Top}_*\rw \ccgp$
is constructed. Further, it is shown in \cite{bcd} that the functors $B$ and $\pcl$
induce an equivalence of categories
\begin{equation}\label{sec1a.eq2}
  \ovl{B}:\frac{\ccgp}{\sim}\;\;\simeq\;\;\hcon :\ovl{\pcl}
\end{equation}
between the localization of $\catu$s with respect to weak equivalences and homotopy
category of connected 3-types.

Given a $\catu$ $\gcl$ it is shown in \cite{bcd} that there is a zig-zag of weak
equivalences in $\ccgp$ between $\gcl$ and $\pcl B\gcl$. We denote by
$[\mi]:\ccgp\rw\ccgp/\!\sim$ the localization functor. Hence $[\gcl]=[\pcl B \gcl]$ for
each $\gcl\in\ccgp$. In general, $[\gcl]=[\gcl']$, if and only if there is a zig-zag of
weak equivalences in $\ccgp$ between $\gcl$ and $\gcl'$.

%%%%%%%%%%%%%%%%%%%%%%%%%%%%%%%%%%%%%%%%%%%%%%%%%%%%%%%

\section{Internal categories in categories of groups with operations.}\label{sec2}
In this section we prove a general fact about internal categories in a category $\ccl$ of
groups with operations. We first recall some preliminary notions.

Internal categories in a category $\ccl$ of groups with operations are internal
groupoids, and are equivalent to reflexive graphs
\begin{equation*}
  C_1\;\pile{\rTo^{d_0}\\ \rTo^{d_1}\\ \lTo_{\sigma_0}}\;C_0
\end{equation*}
satisfying the condition $[\ker d_0,\ker d_1]=1$.

There is a full and faithful nerve functor
\begin{equation*}
  \ncl:\cat\ccl\rTo\cop
\end{equation*}
which has a left adjoint $\pcl:\;\cop\rw\cat\ccl$ given by
\begin{equation*}
  \pcl(H_*):\;H_1/d_2(H_2)\;\pile{\rTo^{\ovl{d}_0}\\ \rTo^{\ovl{d}_1}\\
  \lTo_{\ovl{\sigma}_0}}\;H_0.
\end{equation*}
The unit of the adjunction $u_{H_*}:H_*\rw\ncl\pcl(H_*)$ induces isomorphisms of homotopy
groups $\pi_i$ for $i=0,1$. We say that a map $f$ in $\cat\ccl$ is a weak equivalence if
$\ncl f$ is a weak equivalence in $\cop$; that is, it induces isomorphisms of homotopy
groups. We denote by $\cat \ccl/\!\sim$ the localization of $\cat\ccl$ with respect to
the weak equivalences and by $[\cdot]:\cat\ccl\rw\cat\ccl/\!\sim$ the localization
functor.

The above notion of weak equivalence in $\cop$ is part of the Quillen model category
structure on simplicial objects in an algebraic category given in \cite{qui1},
\cite{qui2}; in this model structure, if $X_*\in\cop$ is cofibrant, then each $X_n$ is
projective in $\ccl$ with respect to the class of regular epimorphisms (see \cite{qui1}).
\begin{lemma}\label{sec2.lem1}
    Let $\ccl$ be a category of groups with operations and let $\gcl\in\cat\ccl$. There
    exists $\phi_{\mathcal{G}}\in\cat\ccl$ whose object of objects $\phi_0$ is projective in $\ccl$ and
    a weak equivalence $u_\gcl:\phi_{\mathcal{G}}\rw\gcl$ in $\cat\ccl$. Further, given a
    morphism $f:\mathcal{G}\rw\mathcal{G'}$ in $\mathrm{Cat}\mathcal{C}$, then there is a
    morphism $\phi_f:\phi_{\mathcal{G}}\rw \phi_\mathcal{G'}$ making the following diagram
    commute
\begin{equation}\label{sec2.eq0}
\begin{diagram}[h=3em, w=2em]
    \phi_{\gcl} & \rTo^{\phi_f} & \phi_{\gcl'}\\
    \dTo^{u_\gcl} && \dTo_{u_{\gcl'}} \\
     \gcl & \rTo^f & \gcl' .
 \end{diagram}
\end{equation}
 Let
$\gcl\supar{f}\gcl'\supar{f'}\gcl''$ be morphisms in $\cat\ccl$. Then
$[\phi_{f'f}]=[\phi_{f'}][\phi_f]$ and $[\phi_\id]=[\id]$.
\end{lemma}
\prf Let $c:\chi\rw\ncl\gcl$ be a cofibrant replacement in $\ccl^{\dop}$. We have a
commutative diagram in $\ccl^{\dop}$
\begin{diagram}[w=4em]
    \ncl\pcl\chi & \rTo^{\ncl\pcl c}& \ncl\pcl\ncl\gcl\\
    \uTo^{u_\chi} && \uTo_{u_{\ncl\gcl}}\\
    \chi & \rTo_c & \ncl\gcl .
\end{diagram}
Since $\pcl\ncl=\id$, $\ncl\pcl\ncl\gcl=\ncl\gcl$ and $u_{\ncl\gcl}=\id$. Since
$\pi_i\ncl\gcl=0=\pi_i\chi$ for $i>1$ and $u_\chi$ induces isomorphisms of $\pi_0,\pi_1$,
it follows that $u_\chi$ is a weak equivalence. Hence, in the diagram above, $c,\;u_\chi$
and $u_{\ncl\gcl}$ are weak equivalences. It follows that $\ncl\pcl c$ is a weak
equivalence. Let $\phi_\mathcal{G}=\pcl\chi$, then $\phi_0=(\pcl\chi)_0=\chi_0$ is
projective in $\ccl$ (since $\chi$ is cofibrant in $\ccl^{\dop}$) and $\pcl
c:\pcl\chi\rw\gcl$ is a weak equivalence in $\cat\ccl$. Let $\nu_f$ be a cofibrant
approximation of the map $\ncl f$, so that the following diagram commutes:
\begin{diagram}[h=3em, w=3em]
     \chi& \rTo^{\nu_f} & \chi'\\
    \dTo^{c} && \dTo_{c'} \\
     \mathcal{NG} & \rTo^{\mathcal{N}f} & \mathcal{NG'} .
\end{diagram}
Applying to this diagram the functor $\mathcal{P}$, (\ref{sec2.eq0}) follows, with
$\phi_f=\pcl\nu_f$, $u_\gcl=\pcl c$, $u_{\gcl'}=\pcl c'$. Let $\ccl^{\dop}_{\leq 1}$ be
the full subcategory of $\ccl^{\dop}$ whose objects $\psi$ are such that $\pi_i\psi=0$
for $i>1$ and let $i:\ccl^{\dop}_{\leq 1}\hookrightarrow\ccl^{\dop}$ be the inclusion.
Then the functor $[\cdot]\pcl i:\ccl^{\dop}_{\leq 1}\rw\cat\ccl/\!\sim$ sends weak
equivalences to isomorphisms. Given morphisms in $\cat\ccl$,
$\gcl\supar{f}\gcl'\supar{f'}\gcl''$, let $\chi\supar{\nu_f}\chi'\supar{\nu_{f'}}\chi''$
be as above. Notice that $\chi,\chi',\chi''\in\ccl^{\dop}_{\leq 1}$ and, by \cite[p.
67-68 ]{qui2} there is a right homotopy between $\nu_{f'f}$ and $\nu_{f'}\nu_{f}$ and
between $\nu_\id$ and $\id$. It follows from a general fact \cite[Lemma 8.3.4]{hirs} that
$[\pcl\nu_{f'f}]=[\pcl\nu_{f'}][\pcl\nu_f]$ and $[\pcl\nu_\id]=\id$, that is
$[\phi_{f'f}]=[\phi_{f'}][\phi_f]$ and $[\phi_\id]=[\id]$. \enpr

%%%%%%%%%%%%%%%%%%%%%%%%%%%%%%%%%%%%%%%%%%%%%%%%%%%%%%%%%%%%%%%%%%%%%%%

\section{The discretization functor.}\label{sec3}
We are going to apply the result of Section \ref{sec2} to the case where
$\ccl=\cat(\gp)$. In this case $\cat\ccl$ is the category of cat$^2$-groups and there is
a classifying space functor $B:\ccgp\mi\mathrm{groups}\rw \mathrm{Top}_*$ obtained by
composition
\begin{equation*}
  \ccgp\rTo^\ncl\gp^{\Delta^{2^{op}}}\rTo^{\ner^*}
  \set^{\Delta^{3^{op}}}\rTo^{diag}\set^{\Delta^{op\;}}\rTo^{|\cdot|}\mathrm{Top_*}.
\end{equation*}
Here $\ncl$ is the multinerve, $\ner^*$ is induced by the nerve functor
$\ner:\gp\rw\set^{\Delta^{op}}$, $diag$ is the multidiagonal and $|\cdot|$ is the
geometric realization.
\begin{lemma}\label{sec3.lem1}
    Let $\ccl=\cat(\gp)$ and let $f:\gcl'\rw\gcl$ be a weak equivalence in $\cat(\ccl)$.
    Then $Bf$ is a weak homotopy equivalence.
\end{lemma}
\prf By hypothesis $\ncl\! f$ induces isomorphisms of simplicial groups
$\pi_i(\ncl\gcl'\!)\!\cong\!\pi_i(\ncl \gcl)\;\;i=0,1$, and therefore isomorphisms of
groups
\begin{equation*}
  \pi_i(\ncl\gcl')_{q_*}\cong\pi_i(\ncl\gcl)_{q_*}
\end{equation*}
for each $q\geq 0,\;i=0,1$. The map of simplicial spaces
$\{B(\ncl\gcl')_{q_*}\}\supar{\ovl{f}}\{B(\ncl\gcl)_{q_*}\}$ is therefore a levelwise
weak equivalence, hence it induces a weak homotopy equivalence of classifying spaces.
Since $B\gcl'$ is weakly homotopy equivalent to the classifying space of
$\{B(\ncl\gcl')_{q_*}\}$, and similarly for $\gcl$, the result follows.\enpr

Given a $\catu$ $\chi$, we denote by $\chi_0\in\cat(\gp)$ the object of objects of the
internal category in $\cat(\gp)$ corresponding to $\chi$. From Lemma \ref{sec2.lem1} and
Lemma \ref{sec3.lem1} we immediately deduce the following proposition.
\begin{proposition}\label{sec3.pro1}
    Let $\gcl\in\ccgp$. There exists $\phi\in\ccgp$ with
    $\phi_0$ projective in $\cgp$ such that $B\phi$ is weakly homotopy equivalent to
    $B\gcl$.
\end{proposition}

Projective objects in $\cat(\gp)$ have been characterized in \cite{ccg}; they have the
form:
\begin{equation*}
  \phi_0:\;F_1\rtimes F_2\;\pile{\rTo^{d_0}\\ \rTo^{d_1}\\ \lTo_{i}}\; F_2
\end{equation*}
where $j:F_1\rw F_2$ is a normal inclusion, $F_1,\,F_2,\,F_2/j(F_1)$ are free groups,
$d_0(x,y)=y$, $d_1(x,y)=j(x)y$, $i(z)=(1,z)$ and $F_2$ acts on $F_1$ by $^{x}y=j(x)y
j(x^{-1})$.

Let $\phi_0^d$ denote the discrete internal category
\begin{equation*}
  \phi_0^d:\;F_2/j(F_1)\;\pile{\rTo^{\id}\\ \rTo^{\id}\\ \lTo_{\id}}\; F_2/j(F_1).
\end{equation*}

If $q_0:F_2\rw F_1/j(F_1)$ is the quotient map and $q_1:F_1\rt F_2\rw F_2/j(F_1)$ is
defined by $q_1(x,y)=q_0(y)$, then (as easily checked) the map
$d=(q_1,q_0):\phi_0\rw\phi^d_0$ is a weak equivalence in $\cat(\gp)$. Since $F_2/j(F_1)$
is free, $q_0$ has a section $q'_0:F_2/j(F_1)\rw F_2$, $q_0q'_0=\id$. If $i:F_2\rw F_1\rt
F_2$ is the inclusion, the map $t:(i q_0',q_0'):\phi^d_0\rw\phi_0$ is also a weak
equivalence in $\cat(\gp)$ and $dt=\id$.

Let $\phi\in\ccgp$ with $\phi_0$ is projective in $\ccl$. Let $\phi_0^d,\;d,\;t$ be as
above and let $\ner:\;\cat(\gp)\rw\gp^{\Delta^{op}}$ be the nerve functor; the multinerve
of $\phi$, as a simplicial object in $\gp^{\Delta^{op}}$ has
\begin{equation*}
  (\ncl\phi)_p=
  \begin{cases}
    \ner \phi_0 & p=0, \\
    \ner\phi_1\tiund{\ner\phi_0}\;\overset{p}{\cdots}\;\tiund{\ner\phi_0}
    \ner\phi_1 & p>0.
  \end{cases}
\end{equation*}
Let $\pt_0,\,\pt_1:\ner\phi_1\rw\ner\phi_0$ and $\sigma_0:\ner\phi_0\rw\ner\phi_1$ be
face and degeneracy maps and let $d:\ner\phi_0\rw\ner\phi_0^d$ and
$t:\ner\phi_0^d\rw\ner\phi_0$ be induced by the weak equivalences $d:\phi_0\rw\phi_0^d$
and $t:\phi_0^d\rw\phi_0$.

We construct a bisimplicial group, which we call the \emph{discrete multinerve of}
$\phi$, denoted $ds\ncl\phi$, as follows. As an object of
$(\gp^{\Delta^{op}})^{\Delta^{op}}$ it is given by
\begin{equation*}
  (ds\,\ncl\phi)_p=
  \begin{cases}
    \ner \phi_0^d & p=0, \\
    \ner\phi_1\tiund{\ner\phi_0}\;\overset{p}{\cdots}\;\tiund{\ner\phi_0}
    \ner\phi_1 & p>0.
  \end{cases}
\end{equation*}
The face and degeneracy operators are $d\pt_0$,
$d\pt_1:(ds\,\ncl\phi)_1\rw(ds\,\ncl\phi)_0$, $\sigma_0
t:(ds\,\ncl\phi)_0\rw(ds\,\ncl\phi)_1 $. All other faces and degeneracy operators are as
in $\ncl\phi$. Since $dt=\id$, it is straightforward to check that $ds\,\ncl\phi$ is a
simplicial object in $\gp^{\Delta^{op}}$.

 Notice that $\ncl\phi$ depends on the choice of
the section $t$. However since $F_2\cong F_1\rt F_2/j(F_1)$, two different choices of
section lead to isomorphic bisimplicial groups.

The following Lemma describes the properties of the discrete multinerve.
\begin{lemma}\label{sec3.lem3}
    Let $\phi\in\ccgp$
with $\phi_0$ projective in $\cat(\gp)$. Then
\begin{itemize}
  \item [i)] $B(ds\,\ncl\phi)=B\phi$.
  \item [ii)] For each $n\geq 2$ the Segal map
\begin{equation*}
  (ds\,\ncl\phi)_n\rTo(ds\,\ncl\phi)_1\tiund{(ds\,\ncl\phi)_0}\overset{n}\cdots\tiund{(ds\,\ncl\phi)_0}
  (ds\,\ncl\phi)_1
\end{equation*}
  is a weak equivalence in $\gp^{\Delta^{op}}$.
\end{itemize}
\end{lemma}

\prf \noindent \\ i)\;The classifying space of any bisimplicial group can be obtained by
composition
\begin{equation*}
\gp^{\Delta^{2^{op}}}\rTo^{diag}
  \gp^{\Delta^{op}}\rTo^{N}\set^{\Delta^{2^{op\;}}}\rTo^{diag}
  \set^{\Delta^{op}}\rTo^{|\cdot|}\mathrm{Top_*},
\end{equation*}
where $diag$ are diagonal functors, $N$ is induced by the nerve functor
$\gp\rw\set^{\Delta^{op}}$ and $|\cdot|$ is the geometric realization.

Let $\psi$ be a bisimplicial group and denote by $\delta^{h}_i,\delta^{v}_i,
\sigma^{h}_i, \sigma^{v}_i$ the horizontal and vertical face and degeneracy operators.
Let $U:\gp\rw\set$ be the forgetful functor. Then
\begin{equation*}
  (Ndiag\;\psi)_{pq}=
  \begin{cases}
    \{\cdot\} & q=0, \\
    U\psi_{pp}\times\overset{q}{\cdots}\times U\psi_{pp} & q>0.
  \end{cases}
\end{equation*}
The horizontal face and degeneracies in $Ndiag\;\psi$ are induced by those in
$diag\;\psi$, while the vertical face and degeneracies in $Ndiag\;\psi$ are the ones
given by the nerve construction.

Let $\chi$ be another bisimplicial group, with face and degeneracies
$\mu^{h}_i,\mu^{v}_i, \nu^{h}_i, \nu^{v}_i$. Suppose that, for all $p>0$ and $q\geq0$,
\begin{equation}\label{lem4.3hyp}
\begin{split}
  & \chi_{pq}=\psi_{pq}\\
  & \delta^{h}_i=\mu^{h}_i:\psi_{p,q+1}\rw\psi_{pq}\\
  & \delta^{v}_i=\mu^{v}_i:\psi_{p+1,q}\rw\psi_{pq}\\
  & \sigma^{h}_i=\nu^{h}_i:\psi_{p,q}\rw\psi_{p,q+1}\\
  & \sigma^{v}_i=\nu^{v}_i:\psi_{p,q}\rw\psi_{p+1,q}.
\end{split}
\end{equation}

 We claim that the bisimplicial sets $Ndiag\;\psi$ and $Ndiag\;\chi$ have
the same diagonal. In fact, by hypothesis (\ref{lem4.3hyp}) $\psi_{nn}=\chi_{nn}$ for
$n\geq1$ so that $(Ndiag\;\psi)_{00}=(Ndiag\;\chi)_{00}$ and
$(Ndiag\;\psi)_{pq}=(Ndiag\;\chi)_{pq}$ for $p>0, q\geq0$. Further, by hypothesis
(\ref{lem4.3hyp}) the face and degeneracy operators
$(diag\;\psi)_{n+1}\rw(diag\;\psi)_{n}$ and $(diag\;\psi)_{n}\rw(diag\;\psi)_{n+1}$
coincide with the respective ones for $diag\;\chi$ for $n>0$. This implies that the face
and degeneracy operators
\begin{equation*}
  \begin{split}
  & (Ndiag\;\psi)_{p+1,q}\rw(Ndiag\;\psi)_{p,q}\\
  & (Ndiag\;\psi)_{p,q+1}\rw(Ndiag\;\psi)_{p,q}\\
  & (Ndiag\;\psi)_{p,q}\rw(Ndiag\;\psi)_{p+1,q}\\
  & (Ndiag\;\psi)_{p,q}\rw(Ndiag\;\psi)_{p,q+1}
  \end{split}
\end{equation*}
for $p>0, q\geq0$ coincide with the respective ones for $Ndiag\;\chi$.

Clearly $(diagNdiag\;\psi)_k=(diagNdiag\;\chi)_k$ for all $k\geq0$. From above, all face
and degeneracy maps in positive dimension of $diagNdiag\;\psi$ coincide with the
respective ones for $diagNdiag\;\chi$. The face maps $U\psi_{11}\rw\{\cdot\}$ are unique
as $\{\cdot\}$ is the terminal object and the degeneracy maps $\{\cdot\}\rw U\psi_{11}$
and $\{\cdot\}\rw U\chi_{11}$ coincide as they both send $\{\cdot\}$ to the unit of the
group $\psi_{11}=\chi_{11}$. In conclusion $diagNdiag\;\psi=diagNdiag\;\chi$. It follows
that $B\psi=B\chi$.

Let $\phi$ be a cat$^{2}$-group satisfying the hypothesis of the lemma,
$\mathcal{N}\phi\in\gp^{\Delta^{2^{op}}}$ be its multinerve and $ds\mathcal{N}\phi$ the
discrete multinerve. By definition of $ds\mathcal{N}\phi$, the two bisimplicial groups
$\mathcal{N}\phi$ and $ds\mathcal{N}\phi$ satisfy hypothesis (\ref{lem4.3hyp}). It
follows from above that $B\mathcal{N}\phi=Bds\mathcal{N}\phi$. By definition
$B\phi=B\mathcal{N}\phi$, hence the result.
\smallskip

 \noindent ii)\; We need to show that for $n\geq 2$,
$$\ner\phi_1\tiund{\ner\phi_0}\;\overset{n}{\cdots}\;\tiund{\ner\phi_0}\ner\phi_1\;\text{ and }
\;\ner\phi_1\tiund{\ner\phi_0^d}\;\overset{n}{\cdots}\;\tiund{\ner\phi_0^d}\ner\phi_1$$
are weakly equivalent simplicial groups.

We proceed by induction on $n$. Since $\pt_0,\pt_1:\ner\phi_1\rw\ner\phi_0$ are
fibrations, the pullback $\ner\phi_1\tiund{\ner\phi_0}\ner\phi_1$ is weakly equivalent to
the homotopy pullback (see \cite{hirs}); since every simplicial group is fibrant this is
weakly equivalent to the homotopy limit of the diagram
$$\ner\phi_1\rTo^{\pt_0}\ner\phi_0\lTo^{\pt_1}\ner\phi_1$$

By the homotopy invariance property of homotopy limits (see \cite{hirs}), since
$$\ner\phi_0\rTo^{d}\ner\phi_0^d$$ is a weak equivalence, then
$$\ner\phi_1\tiund{\ner\phi_0}\ner\phi_1\; \text{ and }\;
\ner\phi_1\tiund{\ner\phi_0^d}\ner\phi_1$$ are weakly equivalent.

Inductively, suppose
$$\ner\phi_1\tiund{\ner\phi_0}\;\overset{n}{\cdots}\;\tiund{\ner\phi_0}\ner\phi_1 \;\text{ and
}\; \ner\phi_1\tiund{\ner\phi_0^d}\;\overset{n}{\cdots}\;
\tiund{\ner\phi_0^d}\ner\phi_1$$ are weakly equivalent.

Notice that
$\ner\phi_1\tiund{\ner\phi_0}\;\overset{n+1}{\cdots}\;\tiund{\ner\phi_0}\ner\phi_1$ is
the pullback of the diagram
$\ner\phi_1\tiund{\ner\phi_0}\;\overset{n}{\cdots}\;\tiund{\ner\phi_0}\ner\phi_1 \rw \ner
\phi_0 \leftarrow \ner\phi_1$. As before, this is weakly equivalent to the homotopy limit
of this diagram; the homotopy invariance property of homotopy limits and the induction
hypothesis then imply that
$$\ner\phi_1\tiund{\ner\phi_0}\;\overset{n+1}{\cdots}\;\tiund{\ner\phi_0}\ner\phi_1 \;\;\text{and}
\;\;\ner\phi_1\tiund{\ner\phi_0^d}\;\overset{n+1}{\cdots}\;\tiund{\ner\phi_0^d}$$
$\ner\phi_1$ are weakly equivalent.\enpr

\begin{definition}\label{sec3.def1}
    The category $\dcl$ of internal 2-nerves is the full subcategory of bisimplicial
    groups $\psi:\Delta^{2^{op}}\rw\gp$ such that:
\begin{itemize}
  \item [i)] Each $\psi_{n*}:\Delta^{op}\rw\gp$ is the nerve of an object of $\cat(\gp)$.
  \item [ii)] $\psi_{0*}:\Delta^{op}\rw\gp$ is constant.
  \item [iii)] The Segal maps $\psi_{n*}\rw\psi_{1*}\tiund{\psi_{0*}}\overset{n}\cdots\tiund{\psi_{0*}}
  \psi_{1*}$ are weak equivalences of simplicial groups.
\end{itemize}
\end{definition}
By Lemma \ref{sec3.lem3} ii) the discrete multinerve $ds\,\ncl\phi$ is an internal
2-nerve. We say that a morphism $f$ in $\dcl$ is a weak equivalence if $Bf$ is a weak
homotopy equivalence.

We aim to define a functor, which we call \emph{discretization}
\begin{equation*}
  \ds:\ccgp/\!\sim\rTo\dcl/\!\sim\;.
\end{equation*}
Let $\cat^2(\gp)_p$ be the full subcategory of $\cat^2(\gp)$ whose objects $\phi$ are
such that $\phi_0$ is projective in $\cat(\gp)$. Let
$[\cdot]:\cat^2(\gp)\rw\cat^2(\gp)/\!\sim$ be the localization functor. From Lemma
\ref{sec2.lem1} there is a functor $S:\cat^2(\gp)\rw\cat^2(\gp)_p/\!\sim$ defined on
objects by $S(\gcl)=[\phi_\gcl]$ and on morphisms by $S(f)=[\phi_f]$. Also, by Lemma
\ref{sec2.lem1}, $S$ sends weak equivalences to isomorphisms. Therefore $S$ induces a
unique functor $\ovl{S}:\cat^2(\gp)/\!\sim\;\rw\cat^2(\gp)_p/\!\sim$ with
$\ovl{S}\circ[\cdot]=S$. On the other hand, the discrete multinerve defines a functor
$ds\,\ncl:\cat^2(\gp)_p\rw\dcl$. On objects, this associates to $\phi$ its discrete
multinerve $ds\,\ncl\phi$. Given a morphism in $\cat^2(\gp)_p$ $F:\phi\rw\phi'$, let
$ds\,\ncl F:ds\,\ncl\phi\rw ds\,\ncl\phi'$ be given by $(ds\,\ncl F)_n=F_n$ for $n\geq
1$, $(ds\,\ncl F)_0=\ovl{F}_0$, where $\ovl{F}_0:\phi^d_0\rw\phi'^{d}_0$ is induced by
$F_0$. Since $dt=\id$, it is easily checked that $\ovl{F}_0 d=d' F_0$, $F_0
t=t'\ovl{F}_0$; this implies that $\ovl{F}_0 d \pt_i=d'\pt'_i F_1$, $i=0,1$ and
$F_1\sigma_0 t=\sigma'_0 t'\ovl{F}_0$, so that $ds\ncl F$ is a morphism in $\dcl$. By
Lemma \ref{sec3.lem3} i), $ds\,\ncl$ preserves weak equivalences, hence it induces a
functor $\ovl{ds\,\ncl}:\cat^2(\gp)_p/\!\!\sim\,\rw\dcl/\!\sim$. Define $\ds$ to be the
composite $\ds=\ovl{ds\,\ncl}\circ\ovl{S}$.
\begin{proposition}\label{sec3.pro2}
    There is a commutative diagram
\begin{diagram}
    \frac{\ccgp}{\sim}&&\rTo^{\ds}&&\frac{\dcl}{\sim}\\
    &\rdTo_B& &\ldTo_B\\
    &&\hcon
\end{diagram}
\end{proposition}
\prf Given $[\gcl]\in\ccgp/\sim$, choose a weak equivalence in $\ccgp$, $\phi\rw\gcl$
with $\phi_0$ projective in $\cat(\gp)$. By definition of $\ds$ and by Lemma
\ref{sec3.lem3} i), we have
$$B\,\ds[\gcl]=B[ds\,\ncl\phi]=[B\,ds\,\ncl
\phi]=[B\phi]=[B\gcl]=B[\gcl].\qquad\qquad\qquad\text{\enpr}$$

%%%%%%%%%%%%%%%%%%%%%%%%%%%%%%%%%%%%%%%

\section{From cat$^2$-groups to Tamsamani's weak 3-groupoids.}\label{sec5}
In this section we connect $\catu$s with a subcategory of \tam weak 3-groupoids.
\begin{definition}\label{sec4.def1}
    Let $\hcl$ be the subcategory of \tam weak 3-groupoids $\psi:\Delta^{3^{op}}\rw\set$
satisfying the additional conditions:
\begin{itemize}
  \item [a)] The constant functor $\psi( [0],\mi,\mi):\Delta^{2^{op}}\rw\set$ takes values in the one-element set.
  \item [b)] For each $m\geq 2$ the Segal maps $\psi(
  [m],\mi,\mi)\rw\psi([1],\mi,\mi)\times\overset{m}
  \cdots \times\psi([1],\mi,\mi)$ are bijections.
\end{itemize}
\end{definition}

Note that objects of $\hcl$ are not strict 3-groupoids because, in general, $\psi(
[m],\!\mi,\!\mi)$ are weak, not strict, 2-nerves. We say a morphism $f:\psi\rw\psi'$ in
$\hcl$ is an external equivalence if it is an external equivalence in $\tcl_3$. By
definition this amounts to requiring that $\psi_{1**}\rw\psi'_{1**}$ is an external
equivalence of weak 2-groupoids and that $T^2\psi\rw T^2\psi'$ is an equivalence of
categories.

Let $\scl$ be as in Section \ref{sec1a}; the following Lemma characterizes external
equivalences in $\scl$.

\begin{lemma}\label{sec5.lem1}
\hspace{1cm}
\begin{itemize}
  \item [a)] Let $f:\phi\rw\phi'$ be a morphism of weak 2-groupoids.
  Then $f$ is an external equivalence if and only if $B f$ is a weak homotopy equivalence.
  \item[b)] Let $g:\psi\rw\psi'$ be a morphism in $\scl$. Then $g$ is an external
  equivalence if and only if $B g$ is a weak homotopy equivalence.
\end{itemize}
\end{lemma}
\prf See Appendix.

Let $\hcl o_\scl(\hcl)$ be the full subcategory of $\sext$ whose objects are in $\hcl$.
Notice that $\hext\subseteq \hcl o_\scl(\hcl)$.

\begin{theorem}\label{sec5.the1}
    There is a functor
\begin{equation*}
  F:\frac{\ccgp}{\sim}\rw\frac{\hcl}{\sim^{ext}}
\end{equation*}
making the following diagram commute
\begin{diagram}
    \frac{\ccgp}{\sim} && \rTo^F && \frac{\hcl}{\sim^{ext}}\\
    &\rdTo_B && \ldTo_B\\
   & &\hcon
\end{diagram}
Further, $F$ induces an equivalence of categories
\begin{equation*}
  \frac{\ccgp}{\sim}\simeq \hcl o_\scl(\hcl).
\end{equation*}
\end{theorem}

\prf Let $N:\dcl\rw\set^{\Delta^{3^{op}}}$ be induced by the nerve functor
$\gp\rw\set^{\Delta^{op}}$ and let $U:\gp\rw\set$ be the forgetful functor. If
$\gcl=\{\psi_{**}\}\in\dcl$ then
\begin{equation*}
  N(\gcl)_{pqr}=
  \begin{cases}
    \{\cdot\} & p=0, \\
    U\psi_{qr}\times\overset{p}{\cdots}\times U\psi_{qr} & p>0.
  \end{cases}
\end{equation*}
We are going to show that $N(\gcl)\in \hcl$. We claim that $N( [1],\mi,\mi)=U\psi_{**}$
is a weak 2-groupoid. In fact, since $\psi\in\dcl$, for each $n\geq 0$, $U\psi_{n*}$ is
the nerve of a groupoid and $U\psi_{0*}$ is a constant simplicial set. For $n\geq 2$, the
map $U\psi_{n*}\rw U\psi_{1*}\tiund{U\psi_{0*}}\cdots\tiund{U\psi_{0*}}U\psi_{1*}$
induces isomorphisms for each $i\geq 0$,
\begin{align*}
  &\pi_i BU\psi_{n*}\cong\pi_i B\psi_{n*}\cong\pi_i B(\psi_{1*}\tiund{\psi_{0*}}\overset{n}\cdots
  \tiund{\psi_{0*}}\psi_{1*})=\\
  & =\pi_i B(U\psi_{1*}\tiund{U\psi_{0*}}\overset{n}\cdots
  \tiund{U\psi_{0*}}U\psi_{1*}).
\end{align*}
Hence the groupoids $U\psi_{n*}$ and
$U\psi_{1*}\tiund{U\psi_{0*}}\overset{n}\cdots\tiund{U\psi_{0*}}U\psi_{1*}$ have weakly
homotopy equivalent classifying spaces. On the other hand, being the underlying
simplicial sets of a simplicial group, they are both fibrant and cofibrant. By a general
theorem of model categories \cite{hov}, they are also simplicially homotopy equivalent,
so the corresponding categories are equivalent. We also have
\begin{equation*}
  TU\psi_{n*}\!\cong\! U\pi_0 \psi_{n*}\!\cong\! U\pi_0(\psi_{1*}\tiund{\psi_{0*}}\overset{n}\cdots
  \tiund{\psi_{0*}}\psi_{1*})\cong U(\pi_0\psi_{1*}\tiund{\pi_0\psi_{0*}}\overset{n}\cdots
  \tiund{\pi_0\psi_{0*}}\pi_0\psi_{1*}).
\end{equation*}
Hence $TU\psi_{n*}$ is the underlying simplicial set of the nerve of a $\cat(\gp)$, hence
it is the nerve of a groupoid. This completes the proof that $U\psi_{n*}=N( [1],\mi,\mi)$
is a weak 2-groupoid. Since, for $p\geq 2$, $N([p],\mi,\mi)=N(
[1],\mi,\mi)\times\overset{p}{\cdots}\times N( [1],\mi,\mi)$, $\,N([p],\mi,\mi)$ is a
weak 2-groupoid for all $p$ and condition b) in the definition of $\hcl$ holds. Clearly
condition a) also holds.

In order to show that $N(\gcl)\in\hcl$ it remains to check that $T^2N(\gcl)$ is a
groupoid. This follows from the fact that
\begin{equation*}
  T^2N(\gcl)[p]=
  \begin{cases}
    \{\cdot\} & p=0, \\
    U \pi_0\pi_0\psi\times\overset{p}{\cdots}\times U \pi_0\pi_0\psi & p>0.
  \end{cases}
\end{equation*}
The functor $N$ induces a functor
\begin{equation*}
  \ovl{N}:\dcl/\sim\rTo\hcl/\sim^{ext}
\end{equation*}
with $\ovl{N}[\gcl]=[N(\gcl)]$. In fact, if $f$ is a weak equivalence in $\dcl$, then
$Nf$ is a morphism in $\hcl$ inducing a weak homotopy equivalence of classifying spaces.
Therefore, by Lemma \ref{sec5.lem1} b), $Nf$ is an external equivalence.

Define $F=\ovl{N}\circ\ds$. Let $[\gcl]\in {\ccgp}/{\sim}$ and choose a weak equivalence
of cat$^2$-groups $\phi\rw\gcl$ with $\phi_0$ projective in $\mathrm{Cat}(\gp)$
(Proposition \ref{sec3.pro1}). Then
\begin{align*}
   & BF[\gcl]=B\ovl{N}\ds[\gcl]=B\ovl{N}[ds\,\ncl\phi]=[BNds\,\ncl\phi]= \\
   & =[B\, ds\,\ncl\phi]=[B\phi]=[B\gcl]=B[\gcl].
\end{align*}
This proves the first part of the theorem.

Let $ R:\hcl o_\scl(\hcl)\rw {\ccgp}/{\sim}$ be given by $R[\chi]=\ovl{\pcl}B\chi$ where
\begin{equation*}
  \ovl{\pcl}:\hcon \rw \frac{\ccgp}{\sim}
\end{equation*}
is induced by the fundamental cat$^2$-group functor $\pcl:Top_*\rw\ccgp$ (see
\cite{bcd}). Let $ i:\hext\hookrightarrow\hcl o_\scl(\hcl)$ be the inclusion. We are
going to show that the pair of functors $(iF,R)$ is an equivalence of categories.

Let $[\gcl]\in{\ccgp}/{\sim}$ and let $\phi\rw\gcl$ be as above. By Lemma \ref{sec3.lem3}
i), $B\phi=Bds\ncl\phi$ so that
\begin{align*}
   & R\, i F[\gcl]=\ovl{\pcl}Bi\ovl{N}\ds[\gcl]=\ovl{\pcl}B\ovl{N}[ds\,\ncl\phi]= \\
   & [\pcl B N\,ds\,\ncl\phi]=[\pcl B\,ds\,\ncl\phi]=[\pcl B\phi].
\end{align*}
On the other hand, by \cite{bcd}, $[\pcl B\phi]=[\phi]$; so that, in conclusion,
\begin{equation*}
  R\, i F[\gcl]=[\phi]=[\gcl].
\end{equation*}
Let $\dsp \pp_3:\hcon\rw\scl$ be the fundamental groupoid functor from \cite{tam}. Recall
that $\pp_3$ sends weak equivalences to external equivalences and, further, for any
$\gcl\in\scl$ there is an external equivalence $\gcl\rw\pp_3 B\gcl$ \cite[Proposition
11.4]{tam}.

Let $[\psi]\in\hcl o_\scl\hcl$. Since $B\psi\simeq B\pcl B\psi$ and there is an external
equivalence in $\scl$, $\psi\rw\pp_3 B\psi$ (see \cite{tam}), we deduce
\begin{equation}\label{sec5.eq3}
  [\psi]=[\pp_3 B\psi]=[\pp_3 B\pcl B\psi]
\end{equation}
in $\sext$. Choose a weak equivalence in $\ccgp$ $\phi\rw\pcl B \psi$ with $\phi_0$
projective in $\cat(\gp)$. Then
\begin{equation}\label{sec5.eq4}
  [\pp_3 B \phi]=[\pp_3 B \pcl B\psi]
\end{equation}
in $\sext$.By Lemma \ref{sec3.lem3} i) we have
\begin{equation}\label{sec5.eq5}
  [Nds\ncl\phi]=[\pp_3Bds\ncl\phi]=[\pp_3B\phi]
\end{equation}
in $\sext$.

In conclusion (\ref{sec5.eq3}), (\ref{sec5.eq4}), (\ref{sec5.eq5}) imply
\begin{equation*}
  [\psi]=[N\, ds\,\ncl\phi]
\end{equation*}
in $\hcl o_\scl(\hcl)$. On the other hand we have
\begin{equation*}
  i\,F R[\psi]=i\,F[\pcl B\psi]= i\,\ovl{N}\ds[\pcl B\psi]=i\,\ovl{N}[ds\,\ncl\phi]=
  [{N}\,ds\,\ncl\phi]
\end{equation*}
so that $i\,F R[\psi]=[\psi]$. This completes the proof that $(iF,R)$ is an equivalence
of categories. \enpr

From the previous theorem, we deduce:
\begin{corollary}\label{sec5.cor1}
    Every object of $\scl$ is equivalent to an object of $\hcl$ through
    a zig-zag of external equivalences.
\end{corollary}
\prf Given $\gcl\in \scl$, by Theorem \ref{sec5.the1}, $[B\pcl B\gcl]=[B F \pcl B \gcl]$
in $\hcon$. Therefore by Theorem \ref{sec1a.the1} $[\pp_3 B\pcl B \gcl]=[\pp_3 B F \pcl
B\gcl]$ in $\sext$. On the other hand, as in the proof of Theorem \ref{sec5.the1}, we
have, in $\sext$:
\begin{equation*}
  [\pp_3 B\pcl B\gcl]=[\pp_3 B \gcl]=[\gcl], \quad [\pp_3 B F \pcl B \gcl]=[F\pcl B
  \gcl].
\end{equation*}
Therefore $[\gcl]=[F\pcl B\gcl]$ in $\sext$, and since $F\pcl B\gcl\in\hcl$, this proves
the result.\enpr

\smallskip

Since the subcategory $\hcl$ of $\scl$ is strictly contained in $\scl$, Corollary
\ref{sec5.cor1} gives a refinement of \tam result (Theorem \ref{sec1a.the1}) showing that
every object of $\scl$, representing a connected 3-type, can be ``semistrictified" to an
object of $\hcl$.

\begin{remark}\rm\label{sec5.rem1}
Let $(\tcl_2,\times)$ be the category of \tam weak 2-groupoids equip-ped with the
cartesian monoidal structure. We observe that $\hcl$ is isomorphic to a full subcategory
of the category $Mon\,(\tcl_2,\times)$ of monoids in $(\tcl_2,\times)$. In fact, by
general theory, using the reduced bar construction, one can show that
$Mon\,(\tcl_2,\times)$ is isomorphic to the category of simplicial objects $\phi$ in
$\tcl_2$ such that $\phi_0$ is trivial and the Segal maps are isomorphisms. On the other
hand, by its definition the category $\hcl$ is the full subcategory of simplicial objects
$\phi$ in $\tcl_2$ with $\phi_0$ trivial and  Segal maps isomoprhisms such that $T\phi$
is a groupoid.
\end{remark}
\section{Another\, semistrictification\, of\, Tamsamani's\, weak\, 3-groupoids\, with\, one\, object}\label{sec5a}
In the previous section we showed that every object of $\scl$ can be ``semistrictified"
to an object of $\hcl$ having the same homotopy type. The goal of this section is to show
that we can perform a different semistrictification from $\scl$ to another subcategory
$\kcl$ of $\scl$ in such a way that the homotopy type is preserved; objects of $\kcl$ are
``semistrict" 3-groupoids, but the directions in which the Segal maps are isomorphisms
rather than equivalences are different from the ones for $\hcl$. In the next section,
both $\hcl$ and $\kcl$ will be used in establishing the connection with Gray groupoids.

We start by defining a strictification functor from \tam weak 2-group-oids to \tam strict
2-groupoids, and by establishing its properties.

Let $\bi:\ncl_2\rw\bic$ and $\nu:2\mi{gpd}\rw\tcl^{st}_2$ be as in Section \ref{sec1a}.
Let ${Bigpd}$ be the category of bigroupoids and their homomorphisms. Then $\bi$
restricts to a functor $\bi:\tcl_2\rw{Bigpd}$. Let $st\,:\bic\rw 2\mi{cat}$ be the
strictification functor described in \cite{gps}. Then $st$ restricts to a functor
$st:{Bigpd}\rw 2\mi{gpd}$.
\begin{definition}\label{sec5a.def1}
    Let $St:\tcl_2\rw\tcl_2^{st}$ be the composite functor
\begin{equation*}
  \tcl_2\supar{Bic}{Bigpd}\supar{st}2\mi{gpd}\supar{\nu}\tcl_2^{st}
\end{equation*}
\end{definition}

\begin{proposition}\label{sec5a.pro1}
    Let $\phi_1,\ldots,\phi_n\in\tcl_2$.
\begin{itemize}
  \item [i)] $St:\tcl_2\rw\tcl_2^{st}$ preserves external equivalences.
  \item [ii)] There is an external equivalence $\widetilde{g}:St(\phi_1\times\cdots\times\phi_n)\rw St\,\phi_1
  \times\cdots\times St\,\phi_n$
  \item [iii)] For each $\phi\in\tcl_2$ there is a functorial zig-zag of weak equivalences
  in $\tcl_2$ between $\phi$ and $St\,\phi$.
\end{itemize}
\end{proposition}
\prf

\noindent i) The functor $Bic$ sends external equivalences to biequivalences; from
\cite{gps} the latter are sent by $st$ to equivalences of 2-categories which, in turn,
are sent to external equivalences by $\nu$. Hence the result.

\noindent ii) It is sufficient to prove the statement for $n=2$, from which the general
case follows easily. Since $\nu$ is an isomorphism and $\bi$ preserves products (this
follows from the explicit description of $\bi$ given in \cite{pao}), it is enough to show
that there is an equivalence of 2-groupoids
\begin{equation*}
  g:st\,(\bi\phi_1\times\bi\phi_2)\rw st\,\bi\phi_1\times st\,\bi\phi_2.
\end{equation*}
For any bicategory $X$, let $\eta_X:X\rw st\, X$ be the biequivalence, natural in $X$
defined in \cite{gps}. There is a commutative diagram
\begin{diagram}[h=2em,labelstyle=\scriptstyle]
    &&&&Bic\,\phi_1\times Bic\,\phi_2\\
    &&&\ldTo(4,4)^{pr_1}& \dTo_{\eta_{Bic\,\phi_1\times Bic\,\phi_2}} &\rdTo(4,4)^{pr_2}\\
    &&&& st(Bic\,\phi_1\times Bic\,\phi_2)  &&\\
    &&&\ldTo^{st\,pr_1} &\dTo^g&\rdTo^{st\,pr_2}&&\\
    Bic\,\phi_1&\rTo_{\eta_{Bic\,\phi_1}}& st\,Bic\,\phi_1 &\lTo_{p_1}&
    st\,Bic\,\phi_1\times st\,Bic\,\phi_2
    &\rTo_{p_2}& st\,Bic\,\phi_2 &\lTo_{\eta_{Bic\,\phi_2}}& Bic\, \phi_2
\end{diagram}
where the map $g$ is uniquely determined by $st\,\pr_1,\,st\,\pr_2$
($p_1,p_2,\pr_1,\pr_2$ are product projections). By universality, it follows that
\begin{equation*}
  g\eta_{\bi\phi_1\times\bi\phi_2}=(\eta_{\bi\phi_1},\eta_{\bi\phi_2}),
\end{equation*}
hence
\begin{equation*}
  B g\, B\eta_{\bi\phi_1\times\bi\phi_2}=B(\eta_{\bi\phi_1},\eta_{\bi\phi_2}).
\end{equation*}
Since $B\eta_{\bi\phi_1\times\bi\phi_2}$ and $B(\eta_{\bi\phi_1},\eta_{\bi\phi_2})$ are
weak homotopy equivalences, so is $B g$. But $B g=B\nu g$, so by Lemma \ref{sec5.lem1},
$\nu g$ is an external equivalence, so $g$ is an equivalence of 2-groupoids, as required.

\noindent iii) By Theorem \ref{sec1a.the2} for each $\phi\in\tcl_2$ there is a map
$\phi\rw N\,\bi\phi$, natural in $\phi$, which is a levelwise equivalence of categories,
hence in particular it is a weak homotopy equivalence. Also, by \cite{gps} there is a
biequivalence $\bi\phi\rw st\,\bi\phi$, natural in $\phi$, which gives rise to a weak
homotopy equivalence $N\,\bi\phi\rw N\,st\,\bi\phi$. Also by \cite{pao}, for each
$\psi\in 2\mi\mathrm{cat}$, $\bi\nu\psi=\psi$ and, by Theorem \ref{sec1a.the2}, we have a
weak homotopy equivalence $St\,\phi=\nu\,st\,\bi\phi\rw N\,\bi\nu\,st\,\bi\phi=
N\,st\,\bi\phi$. In conclusion we obtain a functorial zig-zag of weak equivalences
\begin{equation*}
  \phi\rw N\,\bi\phi\rw N\,st\,\bi\phi\leftarrow St\,\phi
\end{equation*}
\enpr
\begin{definition}\label{sec5a.def2}
    Let $\kcl$ be the subcategory of \tam weak 3-groupoids $\psi:\Delta^{3^{op}}\rw\set$
satisfying the additional conditions:
\begin{itemize}
  \item [a)] The constant functor $\psi([0],\mi,\mi):\Delta^{2^{op}}\rw\set$ takes values
  in the one-element set.
  \item [b)] For each $n$ $\psi([n],\mi,\mi):\Delta^{2^{op}}\rw\set$ is an object of
  $\tcl_2^{st}$, that is for $m\geq 2$ the Segal maps $$\psi([n],[m],\mi)\rw \psi([n],[1],\mi)
  \tiund{\psi([n],[0],\mi)}\overset{m}{\cdots}\tiund{\psi([n],[0],\mi)} \psi([n],[1],\mi)
  $$ are bijections.
\end{itemize}
\end{definition}
We say a morphism in $\kcl$ is an external equivalence if it is an external equivalence
in $\tcl_3$. Let $\hcl o_\scl(\kcl)$ be the full subcategory of $\sext$ whose objects are
in $\kcl$. Notice that $\kext\subseteq \hcl o_\scl(\kcl)$.
\begin{theorem}\label{sec5a.the1}
    There is a functor
\begin{equation*}
  \ovl{St}:\ssext\rw\kkext
\end{equation*}
making the following diagram commute:
\begin{equation}\label{sec5a.eq2}
\xymatrix{ \ssext \ar[r]^{\ovl{St}} \ar[dr]_{B} & \kkext \ar[d]^B\\
& \hcon}
\end{equation}
Further, $F$ induces an equivalence of categories
\begin{equation*}
  \ssext\simeq\hcl o_\scl(\kcl).
\end{equation*}
\end{theorem}
\prf Given an object $\psi:\Delta^{3^{op}}\rw\set$ of $\scl$, let
$\ovl{St}\,\psi:\Delta^{2^{op}}\rw\set$ be given by
\begin{equation*}
  (\ovl{St}\,\psi)_{n**}=St\,\psi_{n**}.
\end{equation*}
We claim that $\ovl{St}\,\psi\in\kcl$. Clearly for each $n$
$(\ovl{St}\,\psi)_{n**}\in\tcl_2^{st}$ and $(\ovl{St}\,\psi)_{0**}=\{\cdot\}$. It remains
to check that the Segal maps
\begin{equation*}
  \nu_n:(\ovl{St}\,\psi)_{n**}\rw (\ovl{St}\,\psi)_{1**}\times\overset{n}\cdots\times
  (\ovl{St}\,\psi)_{1**}
\end{equation*}
are external equivalences in $\tcl_2^{st}$. Since $\psi$ is a weak 3-groupoid, the Segal
maps $\delta_n:\psi_{n**}\rw\psi_{1**}\times\overset{n}\cdots\times \psi_{1**}$ are
external equivalences in $\tcl_2$. On the other hand, by definition of the map
${g}:St\,(\psi_{1**}\times\overset{n}\cdots\times \psi_{1**})\rw
St\,\psi_{1**}\times\overset{n}\cdots\times St\,\psi_{1**}$ it is easily checked that
$\nu_n={g}\circ \,St\,\delta_n$. Since, by Proposition \ref{sec5a.pro1}, ${g}$ and
$St\,\delta_n$ are external equivalences, so is $\nu_n$. This completes the proof that
$\ovl{St}\,\psi\in\kcl$.

If $f:\psi\rw \phi$ is an external equivalence in $\scl$, by definition
$f_1:\psi_{1**}\rw \phi_{1**}$ is an external equivalence in $\tcl_2$ and $Tf:T\psi\rw\
T\phi$ is an equivalence of categories. Hence, by Proposition \ref{sec5a.pro1},
$St\,f_1=(\ovl{St}\,f)_1$ is an external equivalence and, since
$(T\psi)_n=\pi_0\psi_{n**}\cong\pi_0\,St\,\psi_{n**}$, we have $T\psi\cong
T\ovl{St}\,\psi$ and similarly $T\phi\cong T\ovl{St}\,\phi$; therefore $T\ovl{St}\,\psi$
is also an equivalence of categories and in conclusion $\ovl{St}\,f$ is an external
equivalence.

Thus $\ovl{St}$ induces a functor $\dsp\ovl{St}:\sext\rw\kext$. We claim that there is a
weak homotopy equivalence of classifying spaces
\begin{equation}\label{sec5a.eq3}
  B\ovl{St}\,\psi\simeq B\psi.
\end{equation}
In fact, by Proposition \ref{sec5a.pro1} iii), for each $n\geq 0$ there is a functorial
zig-zag of weak equivalences between $\psi_n$ and $St\,\psi_n$. In turn this gives rise
to a zig-zag of maps between $\psi$ and $\ovl{St}\,\psi$; each of these maps is a weak
homotopy equivalence, since it is so levelwise. Thus (\ref{sec5a.eq3}) follows. We
conclude that $[B\ovl{St}\,\psi]=[B\psi]$ in $\dsp\hcon$, showing that (\ref{sec5a.eq2})
commutes.

Let $\kext\overset{i}{\hookrightarrow}\hcl o_\scl(\kcl)\overset{j}{\hookrightarrow}\sext$
be inclusions. We claim that the pair of functors $(i\ovl{St},j)$ is an equivalence of
categories between $\hcl o_\scl(\kcl)$ and $\sext$. In fact, let $[\psi]\in\sext$. Since
$[B\,\ovl{St}\,\psi]=[B\psi]$ in $\dsp\hcon$ we have
\begin{equation*}
  [\pp_3B\ovl{St}\,\psi]=[\pp_3 B\psi]
\end{equation*}
in $\sext$. On the other hand, by \cite{tam}, we have
\begin{equation*}
  [\ovl{St}\,\psi]=[\pp_3B\ovl{St}\,\psi], \quad [\psi]=[\pp_3 B\psi]
\end{equation*}
in $\sext$. Therefore
\begin{equation*}
  j\,i\,\ovl{St}\,[\psi]=[\ovl{St}\,\psi]=[\psi].
\end{equation*}
Finally, let $[\phi]\in\hcl o_\scl(\kcl)$, $\,\phi\in\kcl$. Then
\begin{equation*}
  i\,\ovl{St}\,j[\phi]=i[\ovl{St}\,\phi]=[\ovl{St}\,\phi]=[\phi].
\end{equation*}
Therefore $(i\,\ovl{St}\,,j)$ is an equivalence of categories.\enpr

%%%%%%%%%%%%%%%%%%%%%%%%%%%%%%%%%%%%%%%%%%%%%%%%%%%%%%%%%%%%%%%%%

\section{The comparison with Gray groupoids}\label{sec6}
It is well known that Gray groupoids are semi-strict algebraic models of homotopy
3-types; see \cite{jt}, \cite{le}, \cite{be}, \cite{kp}. Recall that \gra is the category
of 2-categories with monoidal structure given by the Gray tensor product. A Gray category
is a category enriched in \gra. A Gray groupoid is a Gray category whose cells of
positive dimension are invertible. Denote by $Gray\mi gpd_0$ the category of Gray
groupoids with one object.

In this section we compare our semistrict models of connected 3-types to Gray groupoids;
given objects of $\hcl o_\scl(\hcl)$ and of $\hcl o_\scl(\kcl)$, we will associate to
them an object of $\hgra$ representing the same homotopy type.

The functor $\ovl{St}:\scl\rw\kcl$ described in the proof of Theorem \ref{sec5a.the1}
gives by restriction a functor $\ovl{St}:\hcl\rw\kcl$. Let $\ovl{St}\hcl$ be the full
subcategory of $\mathcal{K}$ whose objects have the form $\ovl{St}\phi, \; \phi \in
\mathcal{H}$. Notice that, in general, objects of $\ovl{St}\hcl$ are weak rather than
strict 3-groupoids. In fact, given $\phi\in \hcl$, although the Segal maps
$\phi_{n**}\rw\phi_{1**}\times \overset{n}{\cdots}\times \phi_{1**}$ are bijections, the
maps $St\,(\phi_{n**})\cong St\,(\phi_{1**}\times \overset{n}{\cdots}\times
St\,\phi_{1**})\rw St\,(\phi_{1**})\times \overset{n}{\cdots}\times St\,(\phi_{1**})$ are
merely external equivalences in $\tcl_2^{st}$; this follows from the fact that the
functor $st:{Bicat}\rw2\mi\cat$ does not preserve products strictly but only up to
equivalences of $\twc$, as easily seen from its definition.
\begin{lemma}\label{sec6.lem1}
    Every object of $\kcl$ is equivalent to an object of
    $\ovl{St}\hcl$ through a zig-zag of external equivalences.
\end{lemma}

\prf Let $[\phi]\in\hcl o_\scl(\kcl)$. By \cite{tam}, $[\phi]=[\pp_3B\phi]$ in $\sext$.
Also, $[B\phi]=[B\pcl B\phi]$ in $\dsp\hcon$, hence by \cite{tam} $[\pp_3 B\phi]=[\pp_3
B\pcl B\phi]$ in $\sext$. Therefore
\begin{equation}\label{sec6.eq1}
  [\phi]=[\pp_3 B\pcl B\phi]
\end{equation}
in $\sext$.

By Theorem \ref{sec5.the1}, $[BF\pcl B\phi]=[B\pcl B\phi]$ in $\hcon$, therefore
\begin{equation}\label{sec6.eq2}
  [\pp_3 BF\pcl B\phi]=[\pp_3 B\pcl B\phi]
\end{equation}
in $\sext$.

By Theorem \ref{sec5a.the1}, $[BF\pcl B\phi]=[B\ovl{St}\,F\pcl B\phi]$ in $\hcon$,
therefore
\begin{equation}\label{sec6.eq3}
  [\pp_3 B F \pcl B \phi]=[\pp_3 B \ovl{St}\,F\pcl B \phi]
\end{equation}
in $\sext$. By \cite{tam}
\begin{equation}\label{sec6.eq4}
  [\pp_3 B \ovl{St}\, F \pcl B\phi]=[\ovl{St}\,F\pcl B\phi]
\end{equation}
in $\sext$. Hence (\ref{sec6.eq1}), (\ref{sec6.eq2}), (\ref{sec6.eq3}), (\ref{sec6.eq4})
imply
\begin{equation*}
  [\phi]=[\ovl{St}\,F\pcl B\phi]
\end{equation*}
in $\hcl o_\scl(\kcl)$. Since $\ovl{St}\,F\pcl B\phi\in\ovl{St}\,\hcl$ this proves the
result.\enpr

\begin{theorem}\label{sec6.the1}
    \;There are functors $\;S:\hcl o_\scl(\hcl)\;\rw\;\hgra\;$ and $\;T:\hcl o_\scl(\kcl)\rw\hgra$
    making the following diagram commute:
\begin{equation}\label{sec6.eq5}
\xymatrix{
    \hcl o_\scl(\hcl) \ar[rr]^S \ar[drr]_B&&  \hgra \ar[d]^B  && \hcl o_\scl(\kcl)\ar[ll]_T  \ar[dll]^B\\
    &&  \hcon & }
\end{equation}
\end{theorem}
\prf Let $\phi:\Delta^{3^{op}}\rw\set$ be an object of $\hcl$. Denote
$\phi_n=\phi_{n**}$. By remark \ref{sec5.rem1}, $\hcl$ is isomorphic to a full
subcategory of the category $Mon(\tcl_2,\times)$ of monoids in $(\tcl_2,\times)$. The
monoid corresponding to $\phi$ is $(\phi_1,\pt_{02},\sigma_0)$ where
$\pt_{02}:\phi_2\cong\phi_1\times\phi_1\rw\phi_1$ is the face operator induced by
$[1]\rw[2]$ $0\rw 0$ $1\rw 2$ and $\sigma_0:\phi_0\rw\phi_1$ is a degeneracy operator.

The functor $Bic:\tcl_2\rw {Bicat}$  preserves products and the terminal object, hence it
is a monoidal functor from $(\tcl_2,\times)$ to $({Bicat},\times)$. In \cite{gps} it is
shown that the functor $st\,:{Bicat}\rw 2\mi\mathrm{cat}$ is in fact a monoidal functor
$st\,:({Bicat},\times)\rw \gra$. Denote by
$\eta_{X,Y}:\,st\,X\otimes\,st\,Y\rw\,st\,(X\times Y)$ the tensor product constraint. It
follows that the composite functor
\begin{equation*}
  st\,\circ Bic:(\tcl_2,\times)\rw\gra
\end{equation*}
is monoidal. Therefore, from a general fact, a monoid in $(\tcl_2,\times)$ is sent by
$st\,\circ Bic$ to a monoid in \gra. More precisely, the monoid
$(\phi_1,\pt_{02},\sigma_0)$ in $(\tcl_2,\times)$ is sent to the monoid
$(st\,Bic\,\phi_1,\mu,n)$ in \gra\, where
\begin{equation*}
\mu=(st\,Bic\,\pt_{02})(\eta_{Bic\,\phi_1,Bic\,\phi_1}),\quad n=st\,Bic\,\sigma_0.
\end{equation*}
 The category $Mon(\gra)$ of monoids in \gra is isomorphic to the full subcategory of
Gray categories consisting of those with just one object. Thus we can view
$(st\,Bic\,\phi_1,\mu,n)$ as a Gray category with one object and 2-categorical hom-set
given by $st\,Bic\,\phi_1$; composition is the 2-functor $$\mu:st\,Bic\,\phi_1\otimes
\,st\,Bic\,\phi_1\rw\,st\,Bic\,\phi_1.$$ It is immediate to see that all cells in
positive dimension of this Gray category are invertible, so this is in fact a Gray
groupoid. We denote this Gray groupoid by $\gcl_\phi$.

We are now going to show that there is a weak homotopy equivalence of classifying spaces
$B\gcl_\phi\simeq B\phi$.

Let $W_\phi:\dop\rw 2\mi gpd$ be the reduced bar construction for the monoid
$(st\,\bi\phi_1,\mu,n)$, so that $\nu W_\phi:\dop\rw\tcl_2^{st}$. We observe that there
is a map $l:\nu W_\phi\rw \ovl{St}\,\phi$, where $l_n=\id$ for $n=0,1$ while for $n>1$
$l_n$ is given by
\begin{align*}
   & (\nu W_\phi)_n=\nu(st\,\bi\phi_1\otimes\overset{n}{\cdots}\otimes \,st\,\bi\phi_1) \supar{\nu\eta} \\
   & \nu  st\,(\bi\phi_1\times\overset{n}{\cdots}\times\bi\phi_1)=\nu \,st\,\bi(\phi_1\times
   \overset{n}{\cdots}\times\phi_1)= (\ovl{St}\,\phi)_n.
\end{align*}
We claim that $l$ induces isomorphisms of homotopy groups $\pi_i B\gcl_\phi\cong\pi_i
B\phi$ for all $i\geq 0$. In fact in \cite[p.63]{be} the homotopy groups of a Gray
groupoid $\gcl$ (called there lax 3-groupoid) with respect to a base point $*$ are
described algebraically as follows:
\begin{align*}
   & \pi_1(\gcl)=\mathrm{Aut}_\gcl(*)/\sim \\
   & \pi_2(\gcl)=\mathrm{Aut}_\gcl(1_*)/\sim \\
   & \pi_3(\gcl)=\mathrm{Aut}_\gcl(1_{1_*})
\end{align*}
where the equivalence relation is induced by the cells of the next higher dimension.

When $\gcl=\gcl_\phi$, recalling the algebraic expression of the homotopy groups of a
strict 2-groupoid as in \cite{ms}, we therefore find:
\begin{equation}\label{sec6.eq6}
\begin{split}
   & \pi_1(\gcl_\phi)=\mathrm{Aut}_{\gcl_\phi}(*)/\sim\;\; =\pi_0(st\,Bic\,\phi_1)\\
   & \pi_2(\gcl_\phi)=\mathrm{Aut}_{\gcl_\phi}(1_*)/\sim\;\;\cong\pi_0\hund{st\,Bic\phi_1}(*,*)=
   \pi_1\,st\,Bic\,\phi_1 \\
   & \pi_3(\gcl_\phi)=\mathrm{Aut}_{\gcl_\phi}(1_{1_*})\cong\mathrm{Aut}_{st\,Bic\,\phi_1}(1_*)=
   \pi_2\,st\,Bic\,\phi_1
\end{split}
\end{equation}
By Proposition \ref{sec5a.pro1}, $\pi_i\,st\,Bic\,\phi_1\cong\pi_i\phi_1$ for all $i$
and, by the proof of Lemma \ref{sec5.lem1} and by Proposition \ref{sec5a.pro1},
$\pi_i\,\ovl{St}\,\phi\cong\pi_i\phi\cong\pi_{i-1}\phi_1$ for all $i>0$. Therefore from
(\ref{sec6.eq6}) the claim follows.

We conclude that $[B\gcl_\phi]=[B\phi]$ in $\hcon$. Given ${[\phi]}\in\hcl o_\scl(\hcl)$,
$\phi\in\hcl$, define
\begin{equation*}
  S[\phi]=[\gcl_\phi].
\end{equation*}
This functor is well defined since if $[\phi]=[\phi']$ in $\hcl o_\scl(\hcl)$, then
$\phi$ and $\phi'$ are externally equivalent, hence $B\phi\simeq B\phi'$ so, from above
$B\gcl_\phi\simeq B\gcl_{\phi'}$, which implies $[\gcl_\phi]\simeq [\gcl_{\phi'}]$ in
$\hgra$.

Given $[\psi]\in\hcl o_\scl(\kcl)$, by Lemma \ref{sec6.lem1} there is $\phi\in\hcl$ such
that $[\psi]=[\ovl{St}\,\phi]$ in $\hcl o_\scl(\mathcal{H})$. Define
\begin{equation*}
  T[\psi]=T[\ovl{St}\,\phi]=[\gcl_\phi].
\end{equation*}
This is well defined because if $[\ovl{St}\,\phi]=[\ovl{St}\,\phi']$  then
$B\,\ovl{St}\,\phi=B\,\ovl{St}\,\phi'$ hence $B\phi\simeq B\phi'$ so $B\gcl_\phi\simeq
B\gcl_{\phi'}$, which implies $[\gcl_\phi]\simeq [\gcl_{\phi'}]$. It is immediate that
(\ref{sec6.eq5}) commutes.\enpr
\smallskip

%%%%%%%%%%%%%%%%%%%%%%%%%%%%%%%%%%%%%%%%%%%%%%%%%%%%%%%%%%%%%%%%%

\section*{Appendix }
\noindent\textbf{Proof of Lemma \ref{sec5.lem1}}

a) Let $f$ be an external equivalence. The fact that $B f$ is then a weak homotopy
equivalence is proved in \cite[Proposition 11.2]{tam}. Suppose, conversely, that $B f$ is
a weak homotopy equivalence. Since $\phi$ is a weak 2-groupoid, $\pi_1(\phi)$ is the
nerve of a group and $\pi_n(\phi)=0$ for $n>1$, then $\phi$ satisfies the $\pi_t$-Kan
condition in the sense of \cite{bouf}, for all $t\geq 0$. By \cite[Theorem B.5]{bouf}
there is a first quadrant spectral sequence
\begin{equation*}
  E^2_{s,t}=\pi_s \pi_t\phi\Rightarrow \pi_{s+t} D
\end{equation*}
where $D=\,diag\,\phi$. Since each $\phi_{n*}$ is the nerve of a groupoid,
$\pi^v_t\phi=0$ for $t>1$, so that $E^2_{s,t}=0$ unless $t=0,1$. It follows that there is
a long exact sequence
\begin{equation*}
  \cdots\rTo\pi_{p+1}D\rTo E^2_{p+1,0}\rTo E^2_{p-1,1}\rTo\pi_p D\rTo E^2_{p,0}\rTo\cdots.
\end{equation*}
On the other hand, since $\phi$ is a weak 2-groupoid, for $n>1$
\begin{equation*}
  \pi_0\phi_{n*}\cong\pi_0(\phi_{1*}\tiund{\phi_{0*}}\overset{n}\cdots\tiund{\phi_0*}\phi_{1*})\cong
    \pi_0\phi_{1*}\tiund{\pi_0\phi_{0*}}\overset{n}\cdots\tiund{\pi_0\phi_0*}\pi_0\phi_{1*}
\end{equation*}
\begin{equation*}
  \pi_1\phi_{n*}\cong\pi_1(\phi_{1*}\tiund{\phi_{0*}}\overset{n}\cdots\tiund{\phi_0*}\phi_{1*})\cong
    \pi_1\phi_{1*}\tiund{\pi_1\phi_{1*}}\overset{n}\cdots\tiund{\pi_1\phi_1*}\pi_1\phi_{1*}.
\end{equation*}
It follows that
\begin{align*}
  E^2_{s,0}= & \pi_s\pi_0\phi=0\quad \text{for } s>1 \\
  E^2_{s,0}= & \pi_s\pi_1\phi=0\quad \text{for } s=0 \;\;\text{and } s>1.
\end{align*}
From the above long exact sequence we deduce
\begin{equation}\label{app.eq1}
  \pi_p D=
  \begin{cases}
    E^2_{0,0}=\pi_0\pi_0\phi &\quad p=0, \\
    E^2_{1,0}=\pi_1\pi_0\phi & \quad p=1, \\
    E^2_{1,1}=\pi_1\pi_1\phi=\pi_1\phi_{1*}&\quad  p=2, \\
    0 & \quad p>2.
  \end{cases}
\end{equation}
By hypothesis, $f$ induces isomorphisms of homotopy groups $\pi_i B\phi\cong\pi_i
B\phi'$, $i\geq 0$. Since $\pi_i B\phi=\pi_i\,diag\,\phi$ from above we obtain
isomorphisms
\begin{equation}\label{sec5.eq1}
  \begin{split}
  \pi_0\pi_0\phi\cong & \pi_0\pi_0\,\phi' \\
  \pi_1\pi_0\phi\cong & \pi_1\pi_0\,\phi' \\
  \pi_0\phi_{1*}\cong & \pi_0\,\phi'_{1*}.
  \end{split}
\end{equation}
The first two isomorphisms in (\ref{sec5.eq1})  show that the map
$Tf:\pi_0\phi\rw\pi_0\phi'$ is a weak homotopy equivalence. On the other hand, since
$\pi_0\phi$ and $\pi_0\phi'$ are nerves of groupoids, they are fibrant simplicial sets.
Since every simplicial set is cofibrant, $Tf$ is a weak homotopy equivalence between
fibrant and cofibrant simplicial sets, Thus, by a general result on model categories
\cite{hov}, $Tf$ is a simplicial homotopy equivalence. Hence $Tf$ corresponds to an
equivalence of categories.

In order to show that $f$ is an external equivalence of weak 2-groupoids it remains to
check that, for all $x,y\in\phi_{0*}$ the map
\begin{equation*}
  f_{(x,y)}:\phi_{(x,y)}\rTo\phi'_{(fx,fy)}
\end{equation*}
is an equivalence of categories. Since $\phi_{1*}$ is a groupoid, its first homotopy
group is isomorphic to the endomorphism group of any identity arrow, that is, for any
$z\in\phi_{10}$
\begin{equation*}
  \pi_1\phi_{1*}\cong\hund{\phi_{1*}}(z,z).
\end{equation*}
Further, as $\phi_{10}$ is a groupoid for any $z,z'\in\phi_{10}$ there are isomorphisms
\begin{equation*}
  \hund{\phi_{1*}}(z,z)\cong\hund{\phi_{1*}}(z,z'),
\end{equation*}
and similarly for $\phi'_{1*}$.

Hence by the third isomorphism in (\ref{sec5.eq1}) we deduce that for each
$z,z'\in\phi_{10}$
\begin{equation}\label{sec5.eq2}
   \hund{\phi_{1*}}(z,z')\cong\hund{\phi'_{1*}}(f{z},f{z'}).
\end{equation}
On the other hand recall that
\begin{equation*}
  \underset{x,y\in\phi_{0*}}{\coprod}\phi_{(x,y)}=\phi_{1*},\qquad
  \underset{x',y'\in\phi'_{0*}}{\coprod}\phi_{(x',y')}=\phi'_{1*}.
\end{equation*}
Taking $z,z'$ such that $\pt_0 z=x=\pt_0z'$ and $\pt_1 z=y=\pt_1 z'$ we see that
(\ref{sec5.eq2}) restricts to an isomorphism
\begin{equation*}
   \hund{\phi_{(x,y)}}(z,z')\cong\hund{\phi'_{(fx,fy)}}(f{z},f{z'}).
\end{equation*}
showing that $f_{(x,y)}$ is full and faithful.

Finally let $z\in\phi'_{(fx,fy)}(0)$. Then $fx$ and $fy$ are in the same isomorphism
class of objects of the category $\pi_0\phi'$. Since by (\ref{sec5.eq1})
$\pi_0\pi_0\phi\cong\pi_0\pi_0\phi'$, there exists $w\in\phi_{(x,y)}(0)$ and an
isomorphism $f_{(x,y)}(w)\rw z$. This shows that $f_{(x,y)}$ is essentially surjective on
objects. Hence $f_{(x,y)}$ is an equivalence of categories.

\bigskip

b) The fact that $B g$ is a weak homotopy equivalence when $g$ is an external equivalence
is proved in \cite[Proposition 11.2]{tam}. Suppose, conversely, that $B g$ is a weak
homotopy equivalence. Consider the bisimplicial sets $\phi$ and $\phi'$ defined by
\begin{equation*}
  \phi_{np}= (\,diag\,\psi_{n**})_p \qquad \phi'_{np}=(\,diag\,\psi'_{n**})_p.
\end{equation*}
For each $n\geq 1$ the Segal maps
$\delta_n:\psi_{n**}\rw\psi_{1**}\times\overset{n}{\cdots}\times\psi_{1**}$ induce maps
of simplicial sets
$\ovl{\delta}_n:\phi_{n**}\rw\phi_{1*}\times\overset{n}{\cdots}\times\phi_{1*}$. Since
$\psi,\in\scl$, $\delta_n$ is an external equivalence of weak 2-groupoids, hence
$B\delta_n$ is a weak homotopy equivalence. Since $B\psi_{n**}\simeq B\phi_{n*}$,
$B\ovl{\delta}_n$ is also a weak homotopy equivalence. Therefore, for each $i\geq 0$
\begin{equation*}
  \pi_i\phi_{n*}=
  \begin{cases}
    0 & n=0, \\
    \pi_i\phi_{1*}\times\overset{p}{\cdots}\times\pi_i\phi_{1*} & n>0.
  \end{cases}
\end{equation*}
Thus $\pi_i\phi_{1*}$ is the nerve of a group, so it is fibrant. Hence the bisimplicial
set $\phi$ satisfies the $\pi_*$-Kan condition. By \cite{bouf} there is a spectral
sequence
\begin{equation*}
  E^2_{s,t}=\pi_s\pi_t\phi\Rightarrow\pi_{s+t}(diag\,\phi).
\end{equation*}
We have
\begin{equation*}
  \pi_s\pi_t\phi=
  \begin{cases}
    0 & s=0,\,\,s>1, \\
    \pi_t\phi_{1*} & s=1.
  \end{cases}
\end{equation*}
It follows that
\begin{equation*}
  \pi_n\,diag\,\phi\cong
  \begin{cases}
    0 & n=0, \\
    E^2_{1,n-1}=\pi_{n-1}\phi_{1*} & n>0.
  \end{cases}
\end{equation*}
Similarly we have
\begin{equation*}
  \pi_n\,diag\,\phi'\cong
  \begin{cases}
    0 & n=0, \\
    E^2_{1,n-1}=\pi_{n-1}\phi'_{1*} & n>0.
  \end{cases}
\end{equation*}
By hypothesis, $g$ induces isomorphisms of homotopy groups $\pi_i B\psi\cong\pi_i B
\psi'$ for all $i$. Since $\pi_{i}B\psi\cong\pi_i B\phi\cong\pi_i \,diag\, \phi$ and
similarly for $\psi'$, from above we deduce $\pi_i\phi_{1*}\cong\pi_i\phi'_{1*}$ and
therefore $\pi_i B\psi_{1**}\cong\pi_i B\psi'_{1**}$ for all $i$. Hence the map
$g_|:\psi_{1**}\rw\psi'_{1**}$ of weak 2-groupoids is a weak homotopy equivalence. By a),
it is also an external equivalence of weak 2-groupoids.

To conclude the proof that $g$ is an external equivalence of weak 3-groupoids, it remains
to show that $T^2g$ is an equivalence of categories. We have
\begin{equation*}
  T\psi([p],\mi,\mi)=
  \begin{cases}
    \{\cdot\} & p=0, \\
    \pi_0\pi_0\psi_{1**}\times\overset{p}{\cdots}\times\pi_0\pi_0\psi_{1**} & p>0.
  \end{cases}
\end{equation*}
and similarly for $\psi'$.

Since $g_|:\psi_{1**}\rw\psi'_{1**}$ is an external equivalence of weak 3-groupoids,
$\pi_0g_|$ is an equivalence of categories, therefore
$\pi_0\pi_0\psi_{1**}\cong\pi_0\pi_0\psi'_{1**}$. So $Tg$ is in fact an isomorphism,
hence in particular an equivalence of the corresponding categories.\enpr

%%%%%%%%%%%%%%%%%%%%%%%%%%%%%%%%%%%%%%%%%%%%%%%%%%%%%%%%%%%%%%%%%%%%%%%%%


\begin{thebibliography}{xxxx}

\bibitem[1]{be}C. Berger, Double loop spaces, braided monoidal categories and algebraic
3-type of space, \emph{Contemporary Mathematics}, \textbf{227} (1999), 49-66.
\bibitem[2]{bouf}A. K. Bousfield, E. M. Friedlander, Homotopy theory of $\Gamma$-spaces,
spectra and bisimplicial sets, \emph{Springer Lecture Notes}, no. 658, (1978).
\bibitem[3]{brr}R. Brown, Groupoids and crossed objects in algebraic topology, \emph{Homology,
Homotopy and Applications} \textbf{1} no 1, (1999), 1-78.
\bibitem[4]{bcd}M. Bullejos, A. M. Cegarra, J. Duskin, On cat$^n$-groups and homotopy types,
\emph{Jour. of Pure and Appl. Algebra}, \textbf{86}, (1993), 134-154.
\bibitem[5]{ccg}P. Carrasco, A.N. Cegarra, A.R. Grandjean, (Co)homology of crossed
modules, \emph{Jour. of Pure and Appl. Algebra}, \textbf{168}, (2002), 147-176.
\bibitem[6]{cg}P. Carrasco, A. M. Cegarra, Group-theoretic algebraic models for homotopy
types, \emph{Journal of Pure and Applied Algebra}, \textbf{75}, n. 3 (1991), 195-235.
\bibitem[7]{gps}R. Gordon, A. J. Power, R. Street, \emph{Coherence for tricategories},
Memoirs of the American Mathematical Society, \textbf{117}, no. 558, (1995).
\bibitem[8]{hirs}P. S. Hirschorn, \emph{Model categories and their localizations},
Mathematical Surveys and Monographs, American Mathematical Society, (2002).
\bibitem[9]{hov} M. Hovey, \emph{Model Categories}, Mathematical Surveys and Monographs, Vol. 63, American
Mathematical Society, (1998).
\bibitem[10]{jt}A. Joyal, M. Tierney, Algebraic homotopy types, Handwritten lecture
notes, (1984).
\bibitem[11]{kp}K. H. Kamps and T. Porter, 2-groupoid enrichments in homotopy theory and
algebra, \emph{K-theory}, \textbf{25}, 373-409, (2002).
\bibitem[12]{pao} S.Lack, S.Paoli, 2-nerves for bicategories, Preprint
(2006) arXiv:math.CT/0607271
\bibitem[13]{le}O. Leroy, Sur une notion de 3-cat\'{e}gorie adapte\'{e} \`{a} l'homotopie,
Preprint, AGATA, Univ. Montpellier II, (1994).
\bibitem[14]{lod}J. L. Loday, Spaces with finitely many non-trivial homotopy groups, \emph{J. Pure
and Appl. Algebra} \textbf{24} (1982) 179--202.
\bibitem[15]{ms}I. Moerdijk, J. Svensson, Algebraic classification of equivariant homotopy
2-types, I, \emph{Jour.of Pure and Appl. Algebra} \textbf{89} (1993), 187-216.
\bibitem[16]{orz}G. Orzech, Obstruction theory in algebraic categories I and II, \emph{J.
Pure and Applied Algebra} \textbf{2} (1972) 287-314 and 315-340.
\bibitem[17]{por1}T. Porter, Extensions, crossed modules and internal categories in categories of groups
with operations. \emph{Proc. Edinburgh Math. Soc.} (1987) \textbf{30}, 373-381.
\bibitem[18]{por}T. Porter, N-types of simplicial groups and crossed N-cubes,
\emph{Topology}, \textbf{32}, no. 1, (1993), 5-24.
\bibitem[19]{qui1}D.G. Quillen, \emph{Homotopical algebra}, Lecture Notes in Mathematics,
\textbf{43}, Springer-Verlag (1967).
\bibitem[20]{qui}D.G. Quillen, Spectral sequences of a double semi-simplicial group,
\emph{Topology}, (5), (1966), 155-157.
\bibitem[21]{qui2}D.G. Quillen, On the (co)homology of commutative rings,
\emph{Proc. Symps. Pure Math. Vol. XVII,}, (1968), 65-87.
\bibitem[22]{tam}Z. Tamsamani, Sur des notions de n-cat\'{e}gorie et n-groupoide
non-strictes via des ensembles multi-simpliciaux, \emph{K-theory}, \textbf{16}, (1999),
51-99.
\bibitem[23]{tam1}Z. Tamsamani, Sur des notions de n-cat\'{e}gorie et n-groupoide non
stricted via des ensembles multi-simpliciaux. Thesis, Universit\'{e} Paul Sabatier,
Toulouse, (1966) arXiv:alg-geom/9512006 and arXiv:alg-geom/9607010
\bibitem[24]{whit}J. H. C. Whitehead, Combinatorial homotopy II, \emph{Bull. Amer. Math.
Soc.}, \textbf{55}, (1949), 453-496.
\end{thebibliography}
\end{document}